\documentclass[12pt]{amsart}

\usepackage{amsthm}
\usepackage{epsf}
\usepackage{amssymb}
\usepackage{latexsym}
\usepackage{amsmath}

\newtheorem{theorem}{Theorem}[section] 

\newtheorem{cor}[theorem]{Corollary}
\newtheorem{defn}[theorem]{Definition}

\newtheorem{exmpl}[theorem]{Example}
\newtheorem{prop}[theorem]{Proposition}
\newtheorem{rem}[theorem]{Remark}
\newtheorem{thm}[theorem]{Theorem}

\def\C{{\mathbb C}}
\def\F{{\mathbb F}}

\def\Z{{\mathbb Z}}
\def\T{{\mathbb T}}
\def\TT{{\T \times \T}}
\def\P{{\mathbb P}}
\def\semidirect{\rtimes}

\def \s{\sigma}
\def \CP{\Bbb C \Bbb P}

\def \Aff{\operatorname{Aff}}
\def \len{\operatorname{len}}

\def \Cent{\operatorname{Z}}
\newcommand\Gal[1]{{{#1}_{\operatorname{Gal}}}}
\newcommand\GalAff[1]{{{#1}_{\operatorname{Gal}}^{\operatorname{Aff}}}}

\def\ie{{{i.e.}}}

\def\sub{{\subseteq}}

\newcommand\set[1]{{\{{#1}\}}}
\def\llra{{\,\longrightarrow\,}}
\def\suchthat{{\,:\,\,}}
\def\subjectto{{\,|\,}}

\def\pitil{\tilde{{\Pi}}_1}

\newcommand\sg[1]{{\left<{#1}\right>}}

\def\isom{{\:\cong\:}}
\newcommand\Dref[1]{{Definition \ref{#1}}}
\newcommand\Cref[1]{{Corollary \ref{#1}}}
\newcommand\eq[1]{{(\ref{#1})}}
\newcommand\Eq[1]{{Equation (\ref{#1})}}

\newcommand\FIGURE[3][]{{\bigskip ((Here was supposed to be a figure))\bigskip}}
\newcommand\FIGUREx[4][]{{\bigskip ((Here was supposed to be a figure))\bigskip}}

\long\def\forget#1\forgotten{}

\newcommand\begintable[1][] {{}}

\newcommand\card[1]{{|{#1}|}}

\newcommand\minusset{{-}}
\newcommand\software[1]{{\textsf{#1}}}

\newcommand\ver[2][Changed:]{{}} 
\newcommand\CoxY[2][Y]{{\operatorname{C}_{\operatorname{#1}}(#2)}}
\newcommand\ab{{\operatorname{ab}}}
\def\co{{\,{:}\,}}
\def\ra{{\,\rightarrow\,}}
\newcommand\op[1]{{\operatorname{#1}}}

\newif\ifXY 
\XYtrue     

\newif\ifbigmatrices
\bigmatricestrue 

\ifXY
\usepackage{xy}
\fi

\ifXY
\xyoption{all}
\fi

\begin{document}

\title[The fundamental group of Galois Cover of $\TT$]%
{The fundamental group of Galois cover of the surface $\TT$}
\author[Meirav Amram
 \and Mina Teicher
 \and Uzi Vishne
]{ Meirav Amram
 \and Mina Teicher
 \and Uzi Vishne
}%


\def\BIU{department of mathematics, Bar-Ilan university, Ramat-Gan 52900, Israel}

\address{Meirav Amram, Einstein Institute of Mathematics,
Hebrew university, Jerusalem}
\email{ameirav@math.huji.ac.il}
\address{Mina Teicher, \BIU} \email{teicher@math.biu.ac.il}
\address{Uzi Vishne, \BIU}   \email{vishne@math.biu.ac.il}

\renewcommand{\subjclassname}{%
\textup{2000} Mathematics Subject Classification}


\date{Jun. 1, 2007}

\begin{abstract}
This is the final paper in a series of four, concerning the
surface $\TT$ embedded in $\C\P^8$, where $\T$ is a the one
dimensional torus. In this paper we compute the fundamental group
of the Galois cover of the surface with respect to a generic
projection onto $\C\P^2$, and show that it is nilpotent of class
$3$. This is the first time such a group is presented as the
fundamental group of a Galois cover of a surface.
\end{abstract}

\maketitle

\section{Introduction}\label{intro:sec} 

Studying Galois covers of surfaces proved to be a very useful tool
in understanding the structure of moduli spaces of surfaces of
general type. In 1987 (\cite{MoTe2}), such a construction gave the
first counterexample to the `Bogomolov watershed conjecture' that
surfaces of general type with positive signature have infinite
fundamental groups.

The Galois cover $\Gal{X}$ of an algebraic surface $X$ embedded in
a projective space $\C\P^N$ is the Zariski closure of the fibred
product (with respect to a generic projection to $\C\P^2$) of $n =
\deg(X)$ copies of $X$ where the generalized diagonal is excluded.
The fundamental group of the Galois cover was studied in a series
of papers by Moishezon and Teicher, and other authors. References
can be found in the recent paper \cite{new}. In the current paper
we study $\Gal{(\TT)}$, which is a moduli space of surfaces of
general type.

The first step in the investigation of $\pi_1(\Gal{X})$ is a short
exact sequence (see \cite{MoTe2})
\begin{equation}\label{thisG}
 1 \llra \pi_1(\Gal{X}) \llra \pitil \llra S_n \llra 1 %
\end{equation}
where $\pitil$ is a certain quotient of a fundamental group of the
complement of the branch curve $S$ of the projection $f$. A
presentation of $\pitil$ can be obtained from the braid monodromy
of the curve $S$, and the van Kampen theorem. In principle a
presentation of $\pi_1(\Gal{X})$ can be constructed from that of
$\pitil$ via the Reidmeister-Schreier method \cite{MKS}. However
this is not enough, since the number of generators in $\pitil$ is
multiplied by $\card{S_n} = n!$, and moreover it is usually quite
difficult to identify the group from its presentation.

In this paper we compute $\pi_1(\Gal{X})$ for the surface $X = \T
\times  \T$, where $\T$ is the torus of dimension $1$ (\ie\ an
elliptic curve over $\C$),
%
and the embedding is obtained from a fixed embedding $\T
\hookrightarrow \C\P^2$, and the induced Segre map $X
\hookrightarrow \C\P^8$. We fix a generic projection $f \co \C\P^8
\rightarrow \C\P^2$, and let $S$ be its branch curve.

In \cite{thesis}, \cite{AmTe2} and \cite{AmTe3} we gave a finite
presentation for the group $\pitil$ of \Eq{thisG} for the case of
$X = \T\times\T$. In this paper we use this presentation, the
nature of the embedding of $X$ in $\CP^8$, a degeneration of the
embedding $\T\times\T \hookrightarrow \CP^8$ to a union of $18$
planes (see \cite{AmTe2}), to obtain a concise presentation of
$\pi_1(\Gal{X})$. Moreover, we apply group theoretic techniques
from \cite{Cox} based on ideas from the Reidmeister-Schreier
method \cite{MKS}, and geometric methods which were developed in
\cite{ATV} for a certain quotient of $\pi_1(\Gal{X})$, to compute
$\pi_1(\Gal{X})$.

The main result of this paper is that $\pi_1(\Gal{X})$ is
nilpotent of class $3$ (Theorem \ref{main}). In fact we obtain an
explicit description of $\pi_1(\Gal{X})$ as a quotient of a group
$K$ by a central element $p$, the projective relation, computed in
Equation (\ref{pis}). This group $K$ is given as the kernel of an
explicit epimorphism from a group $K^*$ to $\Z^{10}$, where a
presentation for $K^*$ is given in Corollary \ref{Kstar}.

\section{Overview of new results related to $\TT$}

As before, let $S$ be the branch curve of $X$ under a generic
projection $f \co \CP^8 \ra \CP^2$. It is very difficult to deal
with $\pi_1(\C^2-S)$ directly, so we start by considering $S_0$
where $S_0$ is the branch of a degenerated object. We construct a
degeneration of $X$ into $X_0$, a union of planes, where no three
planes meet in a line. The branch curve of the generic projection
of $X_0$, denoted by $S_0$, is a union of lines.

The surface $X_0$ is composed of $n = \deg(X_0) = 18$ planes with
$27$ intersection lines. The curve $S_0$ is the union of the $27$
intersection lines. It has $9$ intersection points, as depicted in
Figure~\ref{fignumlines}. This diagram by itself is a topological
torus: vertices with the same number are identified. The vertices
are numbered from left to right and bottom up, and the lines are
lexicographically numbered from $1$ to $27$. (Incidentally we
comment that this graph is a $(9,6,2,2)$ strongly regular graph,
namely the valency of each of the $9$ points is $6$, and any two
neighborhoods intersect at $2$ points).

\ifXY \ifbigmatrices
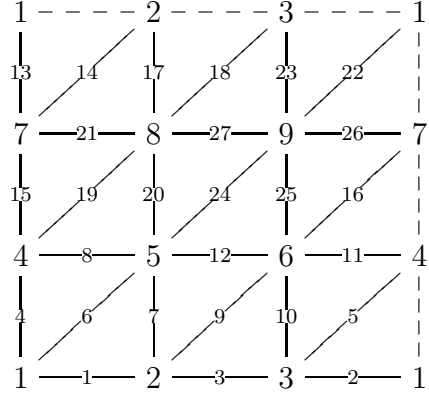
\begin{figure}
\begin{equation*}
\xymatrix@R=30pt@C=36pt{ 
    1 \ar@{--}[r] \ar@{-}[d]|{13}
    & 2 \ar@{--}[r] \ar@{-}[d]|{17} \ar@{-}[ld]|{14}
    & 3 \ar@{--}[r] \ar@{-}[d]|{23} \ar@{-}[ld]|{18}
    & 1 \ar@{--}[d] \ar@{-}[ld]|{22}
\\
    7 \ar@{-}[r]|{21} \ar@{-}[d]|{15}
    & 8 \ar@{-}[r]|{27} \ar@{-}[d]|{20} \ar@{-}[ld]|{19}
    & 9 \ar@{-}[r]|{26} \ar@{-}[d]|{25} \ar@{-}[ld]|{24}
    & 7\ar@{--}[d] \ar@{-}[ld]|{16}
\\
    4 \ar@{-}[r]|{8} \ar@{-}[d]|{4}
    & 5 \ar@{-}[r]|{12} \ar@{-}[d]|{7} \ar@{-}[ld]|{6}
    & 6 \ar@{-}[r]|{11} \ar@{-}[d]|{10} \ar@{-}[ld]|{9}
    & 4 \ar@{--}[d]\ar@{-}[ld]|{5}
\\
    1 \ar@{-}[r]|{1}
    & 2 \ar@{-}[r]|{3}
    & 3 \ar@{-}[r]|{2}
    & 1
}
\end{equation*}
\caption{Numeration of the lines}\label{fignumlines}
\end{figure}

The group $\pi_1(\C^2-S_0)$ is generated by loops
$\Gamma_1,\dots,\Gamma_{27}$, which correspond to the $27$
intersection lines in $S_0$. In the regeneration process
(`inverse' to the degeneration), lines are doubled so $\pi_1(\C^2
- S)$ is generated by $54$ generators, denoted
$\Gamma_1,\Gamma_{1'},\dots,\Gamma_{27},\Gamma_{27'}$.

To compute a complete list of relations, one needs first to
compute the braid monodromy of $S$. To do this, we first compute
the braid monodromy of $S_0$ using the regeneration rules from
\cite{BGT4} and the braid monodromy algorithm for line
arrangements from \cite{BGT2}, Then get a braid monodromy
(factorization) of $S$ from the one of $S_0$ (as in \cite{AmTe2}).
A finite presentation of $\pi_1(\C^2 - S)$ is obtained from the
braid monodromy factorization, via the van~Kampen Theorem
\cite{AmTe3}. A presentation of $\pi_1(\C\P^2-S)$ is obtained by
adding the `projective relation'
\begin{equation}\label{projrel}
\Gamma_{27'} \Gamma_{27} \cdots \Gamma_{1'} \Gamma_{1} = 1.
\end{equation}

The presentation of $\pi_1(\C^2-S)$ involves the above mentioned
$54$ generators, and close to $1700$ relations. This was laid out
in \cite{AmTe3}.
%
A complete presentation is given in the appendix of \cite{thesis},
and can be accessed in [App1] (see the Appendix).

There is a natural projection from $\pi_1(\C\P^2 - S)$ to the
symmetric group $S_{n}$, where $n = 18$. The map is defined by
sending $\Gamma_j$ and $\Gamma_{j'}$ to the transposition $(ab)$
where $a,b$ are the planes intersecting in the line $j$. Since
$(ab)^2 = 1$, the map splits through $\pi_1(\C\P^2 - S)/
\sg{\Gamma_j^2,\Gamma_{j'}^2}$, which we henceforth denote by
$\pitil$. We also let
$$\pitil^{\Aff} = \pi_1(\C^2-S)/ \sg{\Gamma_j^2,\Gamma_{j'}^2}.$$

It is known that $\pi_1(\Gal{X})$ is isomorphic to the kernel of
the projection $\pitil \rightarrow S_{n}$. Our aim in this paper
is to compute the group $\pi_1(\Gal{X})$.

Naturally, it is very difficult to get concrete information on a
group with a presentation of this size. In \cite{ATV} we computed
the quotient $C$ of $\pitil$ obtained by identifying $\Gamma_{j'}
= \Gamma_j$ for every $j$. (The reader might want to check this
paper in order to get a feeling of the type of presentation we are
dealing with, and the difficulties involved in the computation.)

We remark that for most previous surfaces, the quotient $C$
defined in the same manner was computed (implicitly) to be the
symmetric group $S_n$ (for example, see \cite{MoTe2}, \cite{MoTe6}
and \cite{MoRoTe}). The first cases in which $C$ is a larger group
are $\CP^1 \times \T$ (\cite{AGTV} and \cite{thesis}), of the
current surface. See \cite{ATV} for details on this general
approach. In fact these are the only known cases in which the
fundamental group of a Galois cover of a surface is infinite. In
all other cases the fundamental group was shown to be a product of
finite cyclic
groups. In particular $\TT$ is the first case of a group which is not abelian by abelian. %

\medskip



\section{The presentation of $\pitil^{\Aff}$}

We start with a group $\pitil^{\Aff} = \pi_1(\C^2 -
S)/\sg{\Gamma_j^2,\Gamma_{j'}^2}$ generated by the $54$ generators
$\Gamma_1,\Gamma_{1'},\dots,\Gamma_{27},\Gamma_{27'}$, and subject
to about $1700$ relations. These relations are obtained from
applying Van Kapmen's theorem for plane curves to the braid
monodromy factorization of $S$. %
The process is described in details in \cite{AmTe3} and in
\cite{thesis}, and, as mentioned above, the presentation can be
accessed in [App1].
%
However for the sake of this paper it is enough to know that they
all have the following forms:
\begin{enumerate}
\item order two relations: $\Gamma_j^2 = \Gamma_{j'}^2 = 1$,
\item conjugates of the generators commute, such as $[\Gamma_{3},\Gamma_{4}] = 1$ or
$[\Gamma_8,\Gamma_{13}\Gamma_{13'}\Gamma_{14}\Gamma_{13'}\Gamma_{13}]~=~1$,
\item products of conjugates of generators, which have order $3$, such as
$(\Gamma_{15} \cdot \Gamma_{16}\Gamma_{21'}\Gamma_{16})^3 = 1$,
\item $\Gamma_{j'}$ as a conjugate of $\Gamma_j$, 
such as $\Gamma_{9'} =
\Gamma_{3}\Gamma_{3'}\Gamma_{7}\Gamma_{7'}\Gamma_{9}\Gamma_{7'}\Gamma_{7}\Gamma_{3'}\Gamma_{3}$.
\end{enumerate}

As explained above, there is a natural projection from
$\pitil^{\Aff}$ to $S_{18}$. A similar and more abstract situation
was studied in \cite{Cox}. Given a simple connected graph $T$ on
$n$ vertices, a group $\CoxY{T}$ was defined by taking the edges
of $T$ as generators, with the following relations:
\begin{enumerate}
\item[$(R_1)$] Every generator has order $2$ (namely $u^2 = 1$);
\item[$(R_2)$] Two generators commute if their edges are disjoint (so that $[u,v] = 1$);
\item[$(R_3)$] The product of generators $u,v$ has order $3$ if their edges intersect (so $(uv)^3=1$);
\item[$(R_4)$] 
If $u,v,w$ are edges meeting in a point then $[u,vwv] = 1$
(denoting by $u,v,w$ the corresponding generators).
\end{enumerate}
We can rephrase the families $(R_2)$ and $(R_3)$ as $uv = vu$ or
$uvu = vuv$, and in this form they are called the `braid
relations'. It is easy to see that such a group naturally maps
onto $S_n$ (sending a generator $u$ to the transposition $(ab)$
where $a,b$ are the vertices of $u$). The main result of
\cite{Cox} is that $\CoxY{T}$ is a semidirect product of $S_n$ and
a certain normal co-abelian subgroup of $\pi_1(T)^n$, where of
course $\pi_1(T) \isom \F_t$ is the free group on $t$ generators
(\ie\ $t$ is the number of fundamental
cycles in $T$). %
This result was later generalized 
to the case where $T$ is not a simple graph. In this case $(R_3)$
refers to $u,v$ if they intersect in one vertex, $(R_4)$ refers to
three edges $u,v,w$ which intersect in one vertex but cover four
vertices together, and there is a fifth family $(R_5)$ of
relations, associated to every quadruple of edges of which two
edges connect the same two vertices, and each of the other two
edges touches one of these two vertices.

\ifXY
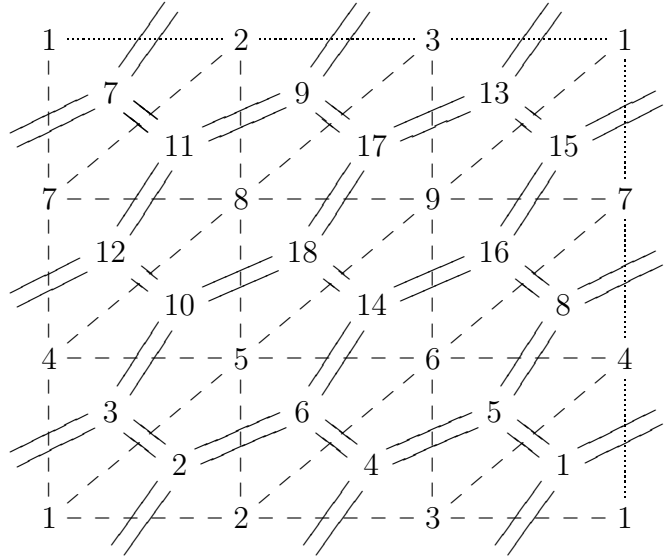
\begin{figure}
\begin{equation}\nonumber
\xymatrix@R=4pt@C=6pt{ 
\\
    {}
    &
    &
    &
    &
    &
    &
    &
    &
    &
    &
    &
    &
{}
\\
    {}
    & 1 \ar@{..}[rrr] \ar@{--}[ddd] 
    &
    &
    & 2 \ar@{..}[rrr] \ar@{--}[ddd] \ar@{--}[lllddd]
    &
    &
    & 3 \ar@{..}[rrr] \ar@{--}[ddd] \ar@{--}[lllddd]
    &
    &
    & 1 \ar@{..}[ddd] \ar@{--}[lllddd]
    &
    &
\\
    {}
    &
    & 7 \ar@{-}@<3pt>[rd] \ar@{-}@<-3pt>[rd] \ar@{-}@<3pt>[rd] \ar@{-}@<-3pt>[rd] \ar@{-}@<3pt>[ruu] \ar@{-}@<-3pt>[ruu] \ar@{-}@<3pt>[lld] \ar@{-}@<-3pt>[lld] 
    &
    &
    & 9 \ar@{-}@<3pt>[rd] \ar@{-}@<-3pt>[rd] \ar@{-}@<-3pt>[ruu] \ar@{-}@<3pt>[ruu] 
    &
    &
    & 13 \ar@{-}@<3pt>[rd] \ar@{-}@<-3pt>[rd] \ar@{-}@<-3pt>[ruu] \ar@{-}@<3pt>[ruu] 
    &
    &
    &
    &
{}
\\
    {}
    &
    &
    & 11 \ar@{-}@<3pt>[rru] \ar@{-}@<-3pt>[rru] \ar@{-}@<3pt>[ldd] \ar@{-}@<-3pt>[ldd]
    &
    &
    & 17 \ar@{-}@<3pt>[rru] \ar@{-}@<-3pt>[rru] \ar@{-}@<3pt>[ldd] \ar@{-}@<-3pt>[ldd]
    &
    &
    & 15 \ar@{-}@<3pt>[rru] \ar@{-}@<-3pt>[rru] \ar@{-}@<3pt>[ldd] \ar@{-}@<-3pt>[ldd]
    &
    &
    &
{}
\\
    {} 
    & 7 \ar@{--}[rrr] \ar@{--}[ddd] 
    &
    &
    & 8 \ar@{--}[rrr] \ar@{--}[ddd] \ar@{--}[lllddd]
    &
    &
    & 9 \ar@{--}[rrr] \ar@{--}[ddd] \ar@{--}[lllddd]
    &
    &
    & 7\ar@{..}[ddd] \ar@{--}[lllddd]
    &
    &
\\
    {}
    &
    & 12 \ar@{-}@<3pt>[rd] \ar@{-}@<-3pt>[rd] \ar@{-}@<3pt>[lld] \ar@{-}@<-3pt>[lld] 
    &
    &
    & 18 \ar@{-}@<3pt>[rd] \ar@{-}@<-3pt>[rd] 
    &
    &
    & 16 \ar@{-}@<3pt>[rd] \ar@{-}@<-3pt>[rd] 
    &
    &
    &
    &
 {}
\\
    {}
    &
    &
    & 10 \ar@{-}@<3pt>[rru] \ar@{-}@<-3pt>[rru] \ar@{-}@<3pt>[ldd]
\ar@{-}@<-3pt>[ldd]
    &
    &
    & 14 \ar@{-}@<3pt>[rru] \ar@{-}@<-3pt>[rru] \ar@{-}@<3pt>[ldd]
\ar@{-}@<-3pt>[ldd]
    &
    &
    & 8 \ar@{-}@<3pt>[rru] \ar@{-}@<-3pt>[rru] \ar@{-}@<3pt>[ldd]
\ar@{-}@<-3pt>[ldd]
    &
    &
    &
{}
\\
    {} 
    & 4 \ar@{--}[rrr] \ar@{--}[ddd] 
    &
    &
    & 5 \ar@{--}[rrr] \ar@{--}[ddd] \ar@{--}[lllddd]
    &
    &
    & 6 \ar@{--}[rrr] \ar@{--}[ddd] \ar@{--}[lllddd]
    &
    &
    & 4\ar@{..}[ddd] \ar@{--}[lllddd]
    &
    &
\\
    {}
    &
    & 3 \ar@{-}@<3pt>[rd] \ar@{-}@<-3pt>[rd] \ar@{-}@<3pt>[lld] \ar@{-}@<-3pt>[lld]
    &
    &
    & 6 \ar@{-}@<3pt>[rd] \ar@{-}@<-3pt>[rd]
    &
    &
    & 5 \ar@{-}@<3pt>[rd] \ar@{-}@<-3pt>[rd]
    &
    &
    &
    &
{}
\\
    {}
    &
    &
    & 2 \ar@{-}@<3pt>[rru] \ar@{-}@<-3pt>[rru] \ar@{-}@<3pt>[ldd] \ar@{-}@<-3pt>[ldd]
    &
    &
    & 4 \ar@{-}@<3pt>[rru] \ar@{-}@<-3pt>[rru] \ar@{-}@<3pt>[ldd] \ar@{-}@<-3pt>[ldd]
    &
    &
    & 1 \ar@{-}@<3pt>[rru] \ar@{-}@<-3pt>[rru] \ar@{-}@<3pt>[ldd] \ar@{-}@<-3pt>[ldd]
    &
    &
    &
{}
\\
    {} 
        & 1 \ar@{--}[rrr] 
    &
    &
    & 2 \ar@{--}[rrr] 
    &
    &
    & 3 \ar@{--}[rrr] 
    &
    &
    & 1 
    &
    &
\\
    {}
    &
    &
    &
    &
    &
    &
    &
    &
    &
    &
    &
    &
{}
}
\end{equation}
\caption{the graph $\hat{T}$} \label{three}
\end{figure}
\else \FIGURE{threereg.ps}{three}
\fi 

Our group $\pitil^{\Aff}$ fits into this general framework, as
follows. The action of $\pitil^{\Aff}$ on the $18$ planes is
defined by $\Gamma_j$ and $\Gamma_{j'}$ switching the two planes
intersecting in the line $j$. Therefore, let $\hat{T}$ be the dual
of the skeleton, namely the graph whose vertices are the planes of
$X_0$, with two edges connecting every two intersecting planes
(one edge for $\Gamma_j$ and one for $\Gamma_{j'}$), depicted as
the solid lines in Figure~\ref{three} (there are $18$ vertices and
$54$ edges). In \cite{ATV} we discovered that $\pitil^{\Aff} /
\sg{\Gamma_{j} = \Gamma_{j'}}$ is a quotient of $\CoxY{T}$ where
$T$ is the graph obtained from $\hat{T}$ by identifying pairs of
edges connecting the same two vertices. Therefore it was expected
that $\pitil^{\Aff}$ is a quotient of $\CoxY{\hat{T}}$.

\subsection{Trivial simplification}

Before we start treating the presentation, we notice that the
computations performed in \cite{thesis} frequently produce
unnecessarily long relations. Taking into account the fact that
all generators have order $2$, we may remove subwords of the form
$\gamma \gamma$ from the relations; furthermore, whenever $\gamma$
and $\delta$ are known to commute, we may replace any subword of
the form $\gamma \delta \gamma$ by $\delta$.

Applying these simple observations reduce the number of distinct
relations from $1692$ (with total length $34644$) to $1599$ (with
total length $18186$).

\subsection{The braid relations}

A possible attack on $\pitil^{\Aff}$ would be to prove that it
satisfies all the five families of relations $(R_1)$--$(R_5)$
defining $\CoxY{\hat{T}}$, making $\pitil^{\Aff}$ a quotient of
$\CoxY{\hat{T}}$, which is very well understood. This would enable
us to get a satisfying description of $\pitil^{\Aff}$. However,
while we do know that $\Gamma_j^2 = \Gamma_{j'}^2 = 1$ and some of
the braid relations are also given, it is difficult to prove the
missing ones (in particular the families $(R_4)$ and $(R_5)$) from
the known relations. We went half-way in this direction, as
follows.

Our first step was to prove that in $\pitil^{\Aff}$ the generators
$\Gamma_j,\Gamma_{j'}$ satisfy the braid relations $(R_2)$ and
$(R_3)$, namely that $\Gamma_{(j)}\Gamma_{(i)}$ has order $2$ or
$3$ depending on whether or not the corresponding edges intersect
(for $i\neq j$). Here, $\Gamma_{(j)}$ stands for either $\Gamma_j$
or $\Gamma_{j'}$. There are $4\binom{27}{2} = 1404$ braid
relations to prove, out of which we expect $4 \cdot 27 \cdot (27 -
5) / 2 = 1188$ pairs to commute (for every $j$ there are four
lines $i$ sharing a common plane with the line $j$); the other $4
\cdot 27 \cdot 4 / 2 = 216$ products are expected to have order
$3$. In practice,
%
%
only $775$ commutators of the expected $1188$ ones, and $152$ of
the expected $216$ products of order $3$, are given as relations
in the presentation obtained above. This leaves $477$ braid
relations which we proved one by one, by hand. This part of the
computation was the most time consuming. The method behind this
lengthy computation is explained in details in \cite{ATV}, where
the same computation was performed on a smaller scale (after
identifying $\Gamma_{j'} =
\Gamma_j$). 

\subsection{Shortening the presentation}\label{ss:short}

Now that the braid relations are known to hold, our next step
involves an automatic simplification of the relations using a
computer program.

It is well known that the word problem is not solvable in a
general presented by generators and relations. Our algorithm is
therefore much less ambitious: given a presentation, we are
looking for relations which can be replaced by shorter ones. For
every two relations, the algorithm is searching for the
intersection of the two words (allowing rotations). Letting $w$
denote the intersection, the two relations can be written as $w
w_1 = w w_2 = 1$ for suitable words $w_1,w_2$. If $\len(w) >
\len(w_1)$, then the relation $w w_2 = 1$ can be replaced by
$w_1^{-1} w_2 = 1$, which is shorter. Furthermore, to somewhat
reduced the chances of getting stuck in local minima, the
algorithm tosses a coin and replaces $w w_2 = 1$ by $w_1^{-1} w_2
= 1$ with probability one half, if $\len(w) = \len(w_1)$.

To speed things up, we use the fact that many pairs of generators
commute. Every few cycles, we scan all the relations for subwords
of the form $\Gamma_j w \Gamma_j^{-1}$, in which all the
generators composing $w$ are known to commute with $\Gamma_j$. In
this case we replace $\Gamma_j w \Gamma_j^{-1}$ by $w$.

The algorithm keeps scanning all pairs of relations as long as
possible reductions are found, and then stops. As an indication of
the performance of this simplistic algorithm in our situation, we
note that the $2331$ relations of total length $21272$ mentioned
before 
were transformed to a list of $1662$
relations, of total length $9312$.

It should be noted that the same algorithm could be ran before
adding the braid relations; however at that stage almost no
shortening occurs since most of the original relations are fairly
complicated.


\subsection{Removing generators}\label{remove}

When the algorithm described above stops, there is no apparent way
to shorten the presentation further. The problem is that in order
to apply what we know on $\CoxY{\hat{T}}$ we still have to prove
the relations in families $(R_4)$ and $(R_5)$, and this is quite
difficult.

Some of the relations express a generator in terms of others, so
it is possible to remove this generator by substitution. For
example, we find that
$$\Gamma_{5'} = \Gamma_{11} \Gamma_{11'} \Gamma_4 \Gamma_{5}
\Gamma_{4'} \Gamma_5 \Gamma_4 \Gamma_{11'} \Gamma_{11}.$$

The `algorithm' applied here was to locate a single generator
$\Gamma_{j'}$ (for $j = 1,\dots,27$) that can be replaced by a
relatively short product, and perform the substitution. Then we
run the shortening algorithm of the previous subsection again. To
save space in what follows, we write $j$ for $\Gamma_j$ and $j'$
for $\Gamma_{j'}$. The first substitution was
\def\sp{{\,\,}}
\def\sps{{\,\,}}
$$27' = 27\sp 18\sp 17'\sp 20\sp 21\sp 19'\sp 21\sp 20\sp 17'\sp
18\sp 27,$$
then
$$26' = 22'\sp 23'\sp 25'\sp 24'\sp 27\sp 24'\sp 25'\sp 23'\sp
22',$$ and then the following:
\def\lineup{-0.25cm}
\begin{eqnarray*}
25' & = & 16'\sp 11'\sp 10'\sp 12'\sp 9'\sp 12'\sp 10'\sp 11'\sp
16'
\\[\lineup]
21' & = & 17\sp 19\sp 19'\sp 21\sp 17'\sp 21\sp 19'\sp 19\sp 17
\\[\lineup]
17' & = & 17\sp 14\sp 1\sp 7\sp 14\sp 17\sp 14\sp 1\sp 7'\sp 1\sp
14\sp 17\sp 14\sp 7\sp 1\sp 14\sp 17
\\[\lineup]
15' & = & 16\sp 21\sp 26\sp 13\sp 14'\sp 13\sp 26\sp 21\sp 16
\\[\lineup]
 9' & = & 3\sp 7\sp 7'\sp 3'\sp 9\sp 3'\sp 7'\sp 7\sp 3
\\[\lineup]
24' & = & 12\sp 7\sp 20\sp 6'\sp 8'\sp 6'\sp 20\sp 7\sp 12
\\[\lineup]
22' & = & 22\sp 2'\sp 13'\sp 13\sp 22\sp 2\sp 22\sp 13\sp 13'\sp
2'\sp 22
\\[\lineup]
20' & = & 8'\sp 6'\sp 6\sp 20\sp 8'\sp 8\sp 8'\sp 20\sp 6\sp 6'\sp
8'\sp
\\[\lineup]
14' & = & 13\sp 21\sp 15\sp 13'\sp 13\sp 16\sp 26\sp 16\sp 13\sp
13'\sp 15\sp 21\sp 13\sp
\\[\lineup]
 7' & = & 6\sp 8'\sp 8\sp 7\sp 6'\sp 7\sp 8\sp 8'\sp 6\sp
\\[\lineup]
 5' & = & 2\sp 10\sp 2'\sp 5\sp 10'\sp 5\sp 2'\sp 10\sp 2\sp
\\[\lineup]
23' & = & 22\sps 13\sps 1\sps 6\sps 4\sps 5\sps 10'\sps 18\sps
17\sps 9\sps 7\sps 14\sps 1\sps 7\sps 9\sps 3 \\[\lineup]
& & \qquad \cdot \ 17\sps 14\sps 17\sps 3\sps 9\sps 7\sps 1\sps
14\sps 7\sps 9\sps 17\sps 18\sps 10'\sps 5\sps 4\sps 6\sps 1\sps
13\sps 22
\\[\lineup]
16' & = & 11\sp 25\sp 24\sp 20\sp 10'\sp 8'\sp 8\sp 8'\sp 6\sp
7\sp 9\sp 7\sp 6\sp 8'\sp 8\sp 8'\sp 10'\sp 20\sp 24\sp 25\sp
11\sp
\\[\lineup]
13' & = & 1\sp 6\sp 4'\sp 13\sp 4\sp 6\sp 1\sp 6\sp 4\sp 13\sp
4'\sp 6\sp 1\sp
\\[\lineup]
12' & = & 7\sp 6\sp 24\sp 20\sp 8'\sp 8\sp 8'\sp 20\sp 24\sp 6\sp
7\sp
\\[\lineup]
 6' & = & 6\sp 1\sp 13\sp 22\sp 4'\sp 6\sp 4\sp 2\sp 6\sp 4'\sp 6\sp 4\sp 2\sp 4\sp 6\sp 22\sp 13\sp 1\sp 6\sp
\\[\lineup]
 2' & = & 22\sp 4\sp 6\sp 13\sp 1\sp 13\sp 6\sp 4\sp 22\sp
\\[\lineup]
 1' & = & 1\sp 6\sp 4\sp 13\sp 2\sp 22\sp 2\sp 13\sp 4\sp 6\sp 1\sp
\\[\lineup]
18' & = & 23\sp 3\sp 10'\sp 2\sp 5\sp 2\sp 10'\sp 3\sp 23\sp
\\[\lineup]
11' & = & 15\sp 19'\sp 8\sp 5\sp 4'\sp 5\sp 8\sp 19'\sp 15\sp
\\[\lineup]
 8' & = & 6\sp 7\sp 12\sp 20\sp 24\sp 20\sp 12\sp 7\sp 6\sp
\\[\lineup]
 3' & = & 17\sp 9\sp 7\sp 14\sp 1\sp 14\sp 7\sp 9\sp 17\sp
\end{eqnarray*}

The original plan was to replace all the generators
$\Gamma_{1'},\dots,\Gamma_{27'}$, so the group will be generated
by $\Gamma_1,\dots,\Gamma_{27}$ which correspond to the simple
graph $T$. However after the above listed substitutions,
$\Gamma_1,\dots,\Gamma_{27}$ remain, together with
$\Gamma_{4'},\Gamma_{10'}$ and $\Gamma_{19'}$. According to the
current relations, each one of these three generators could be
replaced by a product of few hundred generators, a substitution we
did not want to perform.

Fortunately, at this stage there were fairly short substitutions
for $\Gamma_{4}, \Gamma_{10}$ and $\Gamma_{19}$, as follows:
\begin{eqnarray*} 19 & = & 15\sp 8\sp 4\sp 11\sp 5\sp 11\sp 4\sp
8\sp 15\sp
\\[\lineup]
 4 & = & 2\sp 22\sp 6\sp 13\sp 7\sp 9\sp 14\sp 17\sp 3\sp 17\sp 14\sp 9\sp 7\sp 13\sp 6\sp 22\sp 2\sp
\\[\lineup]
10 & = & 5\sp 9\sp 7\sp 1\sp 14\sp 17\sp 18\sp 2\sp 23\sp 2\sp
18\sp 17\sp 14\sp 1\sp 7\sp 9\sp 5. %
\end{eqnarray*}
The group $\pitil^{\Aff}$ is now shown to be generated by
\begin{equation}\label{Delta}
\Delta = \set{\Gamma_j}_{j \neq 4,10,19} \cup
\set{\Gamma_{4'},\Gamma_{10'},\Gamma_{19'}}.
\end{equation}
These generators are in a natural one to one correspondence with
the edges of the simple graph $T$.

\subsection{Quotient of $\CoxY{T}$}\label{ss:quot}

In the previous subsection we expressed half of the original
generators of $\pitil^{\Aff}$ as words on the set $\Delta$,
containing only $27$ generators. During the process, the number of
relations was reduced to $797$ (however the total length went up
to $92610$). The complete presentation on the $27$ generators can
be found in [App2] (see the Appendix). We already know that the
elements of $\Delta$ have order $2$ (so $(R_1)$ is satisfied), and
it is now easy to verify that also satisfy the braid relations
$(R_2)$ and $(R_3)$. In order to prove that $\pitil^{\Aff}$ is a
quotient of $\CoxY{T}$, we only need to show that the fourth
family $(R_4)$ of relations holds (relations from $(R_5)$ no
longer exist, since the current graph of generators, $T$, is
simple).

For example, we need to show that $[\Gamma_{4'},
\Gamma_6\Gamma_8\Gamma_6] = 1$ (lines 4,6, and 8 of the skeleton
bound plane numbered $3$ in Figure~\ref{fignumlines}; dually,
which is the language we prefer here, edges $4,6,8$ of $T$
intersect at the point numbered $3$ in Figure~\ref{three}).
Likewise we need $[\Gamma_1,\Gamma_6\Gamma_7\Gamma_6] = 1$ and so
on, one relation for each of the $18$ triangles of
Figure~\ref{fignumlines}.

It so happens that these relations are all present in the current
list of relations. Thus we proved the claim: $\pitil^{\Aff}$ is a
quotient of $\CoxY{T}$ for the graph $T$ defined before. We note
that $\pi_1(T)$ is the free group of rank $10$, and so $\CoxY{T}$
is a semidirect product of $S_{18}$ and a normal co-abelian
subgroup of $\F_{10}^{18}$ (details are given below).

\newpage
\section{Simplifying the presentation}
\subsection{The structure of $\CoxY{T}$}\label{ss:struc}

\newcommand\lab[1]{#1}

\ifXY
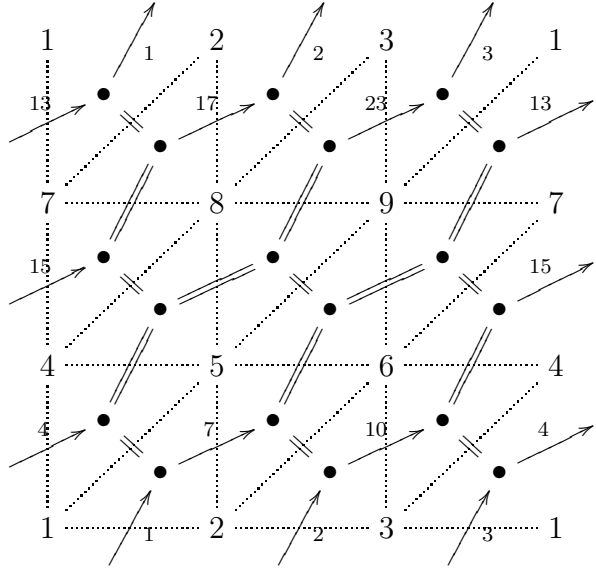
\begin{figure}
\begin{equation}\nonumber
\xymatrix@R=6pt@C=7pt{ 
\\
    {}
    &
    &
    &
    &
    &
    &
    &
    &
    &
    &
    &
    &
{}
\\
    {}
    & 1  \ar@{.}[ddd] 
    &
    &
    & 2  \ar@{.}[ddd] \ar@{.}[lllddd]
    &
    &
    & 3  \ar@{.}[ddd] \ar@{.}[lllddd]
    &
    &
    & 1  \ar@{.}[lllddd]
    &
    &
\\
    {}
    &
    & \bullet \ar@{=}[rd] \ar@{->}[ruu]_{\lab{1}} \ar@{<-}[lld]_{\lab{13}} 
    &
    &
    & \bullet \ar@{=}[rd] \ar@{->}[ruu]_{\lab{2}} 
    &
    &
    & \bullet \ar@{=}[rd] \ar@{->}[ruu]_{\lab{3}} 
    &
    &
    &
    &
{}
\\
    {}
    &
    &
    & \bullet \ar@{->}[rru]^{\lab{17}} \ar@{=}[ldd]
    &
    &
    & \bullet \ar@{->}[rru]^{\lab{23}} \ar@{=}[ldd]
    &
    &
    & \bullet \ar@{->}[rru]^{\lab{13}} \ar@{=}[ldd]
    &
    &
    &
{}
\\
    {} 
    & 7 \ar@{.}[rrr] \ar@{.}[ddd] 
    &
    &
    & 8 \ar@{.}[rrr] \ar@{.}[ddd] \ar@{.}[lllddd]
    &
    &
    & 9 \ar@{.}[rrr] \ar@{.}[ddd] \ar@{.}[lllddd]
    &
    &
    & 7 \ar@{.}[lllddd]
    &
    &
\\
    {}
    &
    & \bullet \ar@{=}[rd] \ar@{<-}[lld]_{\lab{15}} 
    &
    &
    & \bullet \ar@{=}[rd] 
    &
    &
    & \bullet \ar@{=}[rd] 
    &
    &
    &
    &
 {}
\\
    {}
    &
    &
    & \bullet \ar@{=}[rru] \ar@{=}[ldd]
    &
    &
    & \bullet \ar@{=}[rru] \ar@{=}[ldd]
    &
    &
    & \bullet \ar@{->}[rru]^{\lab{15}} \ar@{=}[ldd]
    &
    &
    &
{}
\\
    {} 
    & 4 \ar@{.}[rrr] \ar@{.}[ddd] 
    &
    &
    & 5 \ar@{.}[rrr] \ar@{.}[ddd] \ar@{.}[lllddd]
    &
    &
    & 6 \ar@{.}[rrr] \ar@{.}[ddd] \ar@{.}[lllddd]
    &
    &
    & 4 \ar@{.}[lllddd]
    &
    &
\\
    {}
    &
    & \bullet \ar@{=}[rd] \ar@{<-}[lld]_{\lab{4}}
    &
    &
    & \bullet \ar@{=}[rd]
    &
    &
    & \bullet \ar@{=}[rd]
    &
    &
    &
    &
{}
\\
    {}
    &
    &
    & \bullet \ar@{->}[rru]^{\lab{7}} \ar@{<-}[ldd]^{\lab{1}}
    &
    &
    & \bullet \ar@{->}[rru]^{\lab{10}} \ar@{<-}[ldd]^{\lab{2}}
    &
    &
    & \bullet \ar@{->}[rru]^{\lab{4}} \ar@{<-}[ldd]^{\lab{3}}
    &
    &
    &
{}
\\
    {} 
    & 1 \ar@{.}[rrr] 
    &
    &
    & 2 \ar@{.}[rrr] 
    &
    &
    & 3 \ar@{.}[rrr] 
    &
    &
    & 1 
    &
    &
\\
    {}
    &
    &
    &
    &
    &
    &
    &
    &
    &
    &
    &
    &
{} }
\end{equation}
\caption{Spanning subtree of $T$} \label{subtree}
\end{figure}
\else \FIGURE{threereg.ps}{three}
\fi 

The fundamental group of $T$ is freely generated by $10$
generators. To see this, choose a spanning subtree $T_0$ (which
will contain $18-1 = 17$ edges since $T$ connects $18$ vertices);
then there are $27-17 = 10$ basic cycles, since $T$ has $27$
edges. For the purpose of this paper we choose the subtree to be
the complement of $\Omega = \set{1,2,3,4,7,10,13,15,17,23}$, as in
Figure~\ref{subtree} ($T_0$ is shown in double solid lines).
Although arbitrary, the same choice was made in \cite{ATV}. Figure
\ref{subtree} also provides direction for the edges outside of
$T_0$, to be used later. These arrows are labelled by the edge
numbers, taken from Figure~\ref{fignumlines}.

Recall that we set $n = 18$. For every $\alpha \in
\set{1,\dots,n}$, let $F_\alpha$ denote the free group generated
by the symbols $1_\alpha,2_\alpha,\dots,17_\alpha,23_\alpha$,
namely $\set{\omega_\alpha}_{\omega \in \Omega}$. %
Let $F^*$ denote the direct product of the groups
$F_1,\dots,F_{n}$, so that $F^* \isom (\F_{10})^{n}$. The
symmetric group $S_{n}$ acts on $F^*$ by its action on the
indices.

Let $e_1,e_2,\dots,e_{17},e_{23}$ be generators of a free abelian
group $\Z^{10}$, and define a map $\ab \co F^{*} \rightarrow
\Z^{10}$ by $\ab(\omega_\alpha) = e_{\omega}$ for every $\omega
\in \Omega$ and $\alpha = 1,\dots,n$. Let $F$ denote the kernel of
this map, and note that $\omega_\alpha\omega_\beta^{-1} \in F$ for
every $\omega \in \Omega$ and $\alpha, \beta = 1,\dots,n$; in fact
these elements generate $F$. Obviously $F$ is preserved under the
action of $S_n$.

For an edge $u \in T$, let $\alpha_u,\beta_u$ denote the end
vertices of $u$ (in this order, if $u$ is ordered).
\begin{defn}\label{DefPhi}
Define a map %
\begin{equation}\label{defPhi} %
    \Phi \co \CoxY{T} \ra S_n \semidirect F %
\end{equation} %
by sending $u$ to the transposition $(\alpha_u\beta_u)$ if $u \in
T_0$, and to $(\alpha_\omega \beta_\omega)
\lab{\omega}_{\beta_\omega}^{-1}\lab{\omega}_{\alpha_\omega}$ if
$u = \omega \in \Omega = T \minusset T_0$. %
\end{defn} %

For example, $\Phi(\Gamma_6) = (2 \, 3)$, $\Phi(\Gamma_1) = (2 \,
7) 1_{7}^{-1} 1_{2}$ and $\Phi(\Gamma_7) = (2\, 6) 7_{6}^{-1}
7_2$. In \cite{Cox} it is shown that $\Phi$ is a well defined
isomorphism. Moreover, the natural projection from $\CoxY{T}$ to
$S_{n}$ becomes the projection onto the first component of $S_n
\semidirect F$.

\subsection{$\pitil^{\Aff}$ as a quotient of $S_n \semidirect
F$}\label{FN}

We proved in Subsection \ref{ss:quot} that $\pitil^{\Aff}$ is a
quotient of $\CoxY{T}$, more precisely there is a set $R$ of $797$
elements of $\CoxY{T}$ such that $\pitil^{\Aff}$ is isomorphic to
$\CoxY{T}$ modulo the normal subgroup generated by $R$.

Since the map from $\CoxY{T}$ to $S_n$ splits through
$\pitil^{\Aff}$, $\Phi$ (\Dref{DefPhi}) maps the defining
relations $R$ of $\pitil^{\Aff}$ into $F$. Let $N$ denote the
normal subgroup of $S_n \semidirect F$ generated by these images.
Obviously $N$ is a normal subgroup of $F$, which is moreover
invariant under the action of $S_n$. See Figure~\ref{CYTandPi1}.
It follows that $\Phi$ restricts to an isomorphism
\begin{equation}\label{pitilis} %
\Phi \co \pitil^{\Aff} \ra S_n \semidirect F/N. %
\end{equation}

\ifXY
\begin{figure}
\begin{equation*}
\xymatrix@C=24pt@R=24pt{ 
    F \ar@{^(->}[r]
    \ar@{->}[d]
        & \CoxY{T}        \ar@{->}[r] \ar@{->}[d]
        & S_n \ar@{=}[d]
\\
    F/N \ar@{^(->}[r]
        & \pitil^{\Aff}    \ar@{->}[r]
        & S_n
}
\end{equation*}
\caption{$\CoxY{T}$ and $\pitil^{\Aff}$}\label{CYTandPi1}
\end{figure}
\else
\fi 

We can now rewrite the defining relations in terms of the
generators $\omega_{\alpha}$ of $F$ ($\omega\in \Omega$, $\alpha =
1,\dots,n$). Following the example given after
Definition~\ref{DefPhi}, the relation $(\Gamma_6 \Gamma_7 \Gamma_6
\Gamma_1)^2 = e$ translates to $$((2 \, 3) (2\, 6) 7_{6}^{-1} 7_2
(2 \, 3) (2 \, 7) 1_{7}^{-1} 1_{2})^2 = e,$$ which is equivalent
to $7_{3}^{-1} 7_6 1_{2}^{-1} 1_{7} 7_{6}^{-1} 7_3 1_{7}^{-1}
1_{2} = e$, trivially holding by the definition of $F$.

When this process is applied to the $797$ defining relations,
almost all vanish (namely they express the trivial word in $F$).
We are left with $102$ relations, from which the relation below is
one of the shortest:
$$1_{10}^{-1} 13_{10}^{-1} 2_{10} 4_{10}^{-1} 7_{10} 4_{10}
2_{10}^{-1} 13_{10} 1_{10} 7_{10}^{-1} 7_{14} 1_{14}^{-1}
13_{14}^{-1} 2_{14} 4_{14}^{-1} 7_{14}^{-1} 4_{14} 2_{14}^{-1}
13_{14} 1_{14} = e.$$ 

These $102$ elements generate $N$ (as a normal and invariant
subgroup). Since $N$ is invariant under the action of $S_n$, the
indices in the relations can be made arbitrary. Taking the indices
to be the `generic' $i,j,k,l$ and removing duplicates, the $102$
relations collapse into a list of $66$ relations. This list can be
found in [App3] (see the Appendix).

Next, if a relation has the form $\omega_i u_i \omega_i^{-1} w_j =
e$ for $u_i,\omega_i \in F_i$, $w_j \in F_j$ and $i \neq j$, then
$\omega_i$ commutes with $w_j$ and so it is equivalent to $u_i w_j
= e$. For example,
$$7_i 1_i^{-1} 13_i^{-1} 2_i 4_i^{-1} 7_i 4_i 2_i^{-1} 13_i 1_i 7_i^{-1} 7_i^{-1}
\cdot 7_j 1_j^{-1} 13_j^{-1} 2_j 4_j^{-1} 7_j^{-1} 4_j 2_j^{-1}
13_j 1_j = e$$ becomes %
\begin{equation}\label{rel} %
1_i^{-1} 13_i^{-1} 2_i 4_i^{-1} 7_i 4_i 2_i^{-1} 13_i 1_i 7_i^{-1}
\cdot 7_j 1_j^{-1} 13_j^{-1} 2_j 4_j^{-1} 7_j^{-1} 4_j 2_j^{-1}
13_j 1_j = e. %
\end{equation}
`Cleaning' the $66$ relations from such conjugations (which can be
done by hand), we obtain a list of $49$ relations, all of the
forms $u_i = e$ (a single relation), $u_i \cdot w_j = e$ ($29$
relations), $u_i \cdot w_j \cdot v_k = e$ ($13$ relations) or $u_i
\cdot w_j \cdot v_k \cdot x_l = e$ ($6$ relations) for distinct
$i,j,k,l$ and various elements $u_i \in F_i$, $w_j \in F_j$, $v_k
\in F_k$, $x_{l} \in F_l$.

\forget
\begin{eqnarray*}
17_i^{-1} 1_i 3_i^{-1} 17_i^{-1} 13_i^{-1} 15_i^{-1} 1_i
7_i^{-1} 3_i^{-1} \cdot
 17_j 3_j 1_j^{-1} 17_j 3_j 7_j 1_j^{-1} 15_j 13_j
\\
4_i 2_i^{-1} 13_i 1_i 7_i 1_i^{-1} 13_i^{-1} 15_i^{-1} 17_i^{-1}
1_i \cdot
 4_j^{-1} 1_j^{-1} 17_j 15_j 13_j 1_j 7_j^{-1} 1_j^{-1} 13_j^{-1}
 2_j
\\
4_i^{-1} 7_i 4_i 2_i^{-1} 13_i 1_i 7_i^{-1} 1_i^{-1} 13_i^{-1} 2_i
\cdot
 4_j^{-1} 7_j^{-1} 4_j 2_j^{-1} 13_j 1_j 7_j 1_j^{-1} 13_j^{-1}
 2_j
\\
4_i 2_i^{-1} 13_i 17_i 3_i 7_i^{-1} 3_i^{-1} 17_i^{-1} 13_i^{-1}
2_i 15_i^{-1} \cdot
 4_j^{-1} 15_j 2_j^{-1} 13_j 17_j 3_j 7_j 3_j^{-1} 17_j^{-1} 13_j^{-1} 2_j
\\
 23_i 13_i 17_i 3_i 7_i 1_i^{-1} 13_i^{-1} 23_i^{-1} 17_i^{-1} 1_i 7_i^{-1} 3_i^{-1} \cdot
 23_j 13_j 1_j 7_j^{-1} 3_j^{-1} 17_j^{-1} 13_j^{-1} 23_j^{-1} 3_j 7_j 1_j^{-1} 17_j
\\
 4_i^{-1} 7_i 4_i 2_i^{-1} 13_i 17_i 3_i 7_i^{-1} 3_i^{-1} 17_i^{-1} 13_i^{-1} 2_i \cdot
 17_j 3_j 7_j 3_j^{-1} 17_j^{-1} 13_j^{-1} 2_j 15_j^{-1} 7_j^{-1} 15_j 2_j^{-1} 13_j
\\
 4_i 2_i^{-1} 13_i 1_i 7_i 1_i^{-1} 7_i^{-1} 3_i^{-1} 17_i^{-1} 13_i^{-1} 2_i 13_i 1_i \cdot
 4_j^{-1} 1_j^{-1} 13_j^{-1} 2_j^{-1} 13_j 17_j 3_j 7_j 1_j 7_j^{-1} 1_j^{-1} 13_j^{-1} 2_j
\\
 10_i^{-1} 2_i^{-1} 13_i 17_i 3_i 7_i 1_i^{-1} 13_i^{-1} 2_i 10_i 3_i^{-1} 17_i^{-1} 1_i 7_i^{-1} \cdot
 10_j^{-1} 2_j^{-1} 13_j 1_j 7_j^{-1} 3_j^{-1} 17_j^{-1} 13_j^{-1} 2_j 10_j 7_j 1_j^{-1} 17_j 3_j
\\
 10_i^{-1} 2_i^{-1} 23_i^{-1} 17_i^{-1} 1_i 7_i^{-1} 3_i^{-1} 17_i^{-1} 1_i 4_i 2_i^{-1} 13_i 1_i 7_i \cdot
 10_j 7_j^{-1} 1_j^{-1} 13_j^{-1} 2_j 4_j^{-1} 1_j^{-1} 17_j 3_j 7_j 1_j^{-1} 17_j 23_j 2_j
\\
 4_i 2_i^{-1} 13_i 1_i 7_i^{-1} 3_i^{-1} 17_i^{-1} 13_i^{-1} 2_i 4_i^{-1} 1_i^{-1} 17_i 3_i 7_i \cdot
 4_j^{-1} 7_j^{-1} 3_j^{-1} 17_j^{-1} 1_j 4_j 2_j^{-1} 13_j 17_j 3_j 7_j 1_j^{-1} 13_j^{-1} 2_j
\\
 10_i^{-1} 15_i^{-1} 4_i 10_i 3_i^{-1} 17_i^{-1} 1_i 4_i 2_i^{-1} 13_i 1_i 10_i^{-1} 4_i^{-1} 15_i 7_i \cdot
 10_j 1_j^{-1} 13_j^{-1} 2_j 4_j^{-1} 1_j^{-1} 17_j 3_j 10_j^{-1} 4_j^{-1} 15_j 10_j 7_j^{-1} 15_j^{-1} 4_j
\\
 4_i 2_i^{-1} 13_i 1_i 7_i 1_i^{-1} 15_i 13_i 7_i^{-1} 3_i^{-1} 17_i^{-1} 13_i^{-1} 2_i 15_i^{-1} 1_i \cdot
 4_j^{-1} 1_j^{-1} 15_j 2_j^{-1} 13_j 17_j 3_j 7_j 13_j^{-1} 15_j^{-1} 1_j 7_j^{-1} 1_j^{-1} 13_j^{-1} 2_j
\\
 4_j 2_j^{-1} 13_j 1_j 3_j^{-1} 17_j^{-1} 1_j 7_j^{-1} 1_j^{-1} 13_j^{-1} 2_j 4_j^{-1} 1_j^{-1} 17_j 3_j 7_j \cdot
 4_i 2_i^{-1} 13_i 1_i 7_i 1_i^{-1} 17_i 3_i 1_i^{-1} 13_i^{-1} 2_i 4_i^{-1} 7_i^{-1} 3_i^{-1} 17_i^{-1} 1_i
\\
 10_i^{-1} 15_i^{-1} 4_i 10_i 3_i^{-1} 17_i^{-1} 1_i 4_i 2_i^{-1} 13_i 1_i 7_i^{-1} 10_i^{-1} 7_i \cdot
 10_j^{-1} 4_j^{-1} 15_j 10_j 7_j^{-1} 15_j^{-1} 4_j 10_j 1_j^{-1} 13_j^{-1} 2_j 4_j^{-1} 1_j^{-1} 17_j 3_j \cdot
 10_k 7_k 10_k^{-1} 4_k^{-1} 15_k
\\
 4_i 2_i^{-1} 13_i 17_i 3_i 7_i 1_i^{-1} 13_i^{-1} 15_i^{-1} 7_i^{-1} 7_i^{-1} 3_i^{-1} 17_i^{-1} 13_i^{-1} 2_i 13_i 1_i \cdot
 4_j^{-1} 1_j^{-1} 13_j^{-1} 2_j^{-1} 13_j 17_j 3_j 7_j 7_j 15_j 13_j 1_j 7_j^{-1} 3_j^{-1} 17_j^{-1} 13_j^{-1} 2_j
\\
 4_i 2_i^{-1} 13_i 17_i 3_i 7_i 1_i^{-1} 17_i 3_i 1_i^{-1} 13_i^{-1} 2_i 4_i^{-1} 1_i^{-1} 15_i 13_i 1_i \cdot
 4_j 2_j^{-1} 13_j 1_j 3_j^{-1} 17_j^{-1} 1_j 7_j^{-1} 3_j^{-1} 17_j^{-1} 13_j^{-1} 2_j 4_j^{-1} 1_j^{-1} 13_j^{-1} 15_j^{-1} 1_j
\\
 4_j 2_j^{-1} 13_j 1_j 7_j 4_j 2_j^{-1} 13_j 1_j 7_j 1_j^{-1} 7_j^{-1} 3_j^{-1} 17_j^{-1} 13_j^{-1} 2_j 15_j^{-1} 17_j^{-1} 1_j \cdot
 4_i^{-1} 7_i^{-1} 1_i^{-1} 13_i^{-1} 2_i 4_i^{-1} 1_i^{-1} 17_i 15_i 2_i^{-1} 13_i 17_i 3_i 7_i 1_i 7_i^{-1} 1_i^{-1} 13_i^{-1} 2_i
\\
 4_i 2_i^{-1} 13_i 17_i 3_i 7_i 1_i^{-1} 13_i^{-1} 15_i^{-1} 1_i 7_i^{-1} 3_i^{-1} 17_i^{-1} 13_i^{-1} 2_i 4_i^{-1} 1_i^{-1} 15_i 13_i 1_i \cdot
 4_j 2_j^{-1} 13_j 17_j 3_j 7_j 1_j^{-1} 15_j 13_j 1_j 7_j^{-1} 3_j^{-1} 17_j^{-1} 13_j^{-1} 2_j 4_j^{-1} 1_j^{-1} 13_j^{-1} 15_j^{-1} 1_j
\\
 10_i^{-1} 2_i^{-1} 23_i^{-1} 17_i^{-1} 1_i 4_i 2_i^{-1} 13_i 1_i 7_i 1_i^{-1} 17_i^{-1} 1_i 4_i 2_i^{-1} 13_i 1_i 7_i 1_i^{-1} 17_i 3_i \cdot
 10_j 3_j^{-1} 17_j^{-1} 1_j 7_j^{-1} 1_j^{-1} 13_j^{-1} 2_j 4_j^{-1} 1_j^{-1} 17_j 1_j 7_j^{-1} 1_j^{-1} 13_j^{-1} 2_j 4_j^{-1} 1_j^{-1} 17_j 23_j 2_j
\\
 4_i 2_i^{-1} 13_i 1_i 7_i 1_i^{-1} 7_i^{-1} 3_i^{-1} 17_i^{-1} 13_i^{-1} 2_i 15_i^{-1} 17_i^{-1} 1_i 4_i 2_i^{-1} 13_i 1_i 7_i \cdot
 4_j^{-1} 1_j^{-1} 13_j^{-1} 2_j^{-1} 13_j 17_j 3_j 7_j 1_j 7_j^{-1} 1_j^{-1} 13_j^{-1} 2_j \cdot
 4_k^{-1} 1_k^{-1} 17_k 15_k 13_k 1_k 7_k^{-1} 1_k^{-1} 13_k^{-1} 2_k
\\
 4_j 2_j^{-1} 13_j 1_j 7_j 1_j^{-1} 17_j^{-1} 1_j 4_j 2_j^{-1} 13_j 1_j 7_j 1_j^{-1} 7_j^{-1} 3_j^{-1} 17_j^{-1} 13_j^{-1} 2_j 15_j^{-1} 1_j \cdot
 4_i^{-1} 1_i^{-1} 17_i 1_i 7_i^{-1} 1_i^{-1} 13_i^{-1} 2_i 4_i^{-1} 1_i^{-1} 15_i 2_i^{-1} 13_i 17_i 3_i 7_i 1_i 7_i^{-1} 1_i^{-1} 13_i^{-1} 2_i
\\
 4_i^{-1} 1_i^{-1} 17_i 7_i 7_i 4_i 2_i^{-1} 13_i 1_i 7_i^{-1} 17_i^{-1} 1_i 7_i^{-1} 1_i^{-1} 13_i^{-1} 2_i \cdot
 4_j 2_j^{-1} 13_j 1_j 7_j 1_j^{-1} 7_j^{-1} 3_j^{-1} 17_j^{-1} 13_j^{-1} 2_j 15_j^{-1} 7_j 1_j^{-1} 13_j^{-1} 2_j 4_j^{-1} 7_j^{-1} 7_j^{-1} 15_j 2_j^{-1} 13_j 17_j 3_j 7_j 1_j
\\
 23_i^{-1} 17_i^{-1} 1_i 7_i^{-1} 3_i^{-1} 23_i 13_i 17_i 3_i 7_i 1_i^{-1} 13_i^{-1} \cdot
 23_j 13_j 17_j 3_j 7_j 1_j^{-1} 13_j^{-1} 23_j^{-1} 17_j^{-1} 1_j 7_j^{-1} 3_j^{-1} \cdot
 23_k^{-1} 3_k 7_k 1_k^{-1} 17_k 23_k 13_k 1_k 7_k^{-1} 3_k^{-1} 17_k^{-1} 13_k^{-1} \cdot
 23_l 13_l 1_l 7_l^{-1} 3_l^{-1} 17_l^{-1} 13_l^{-1} 23_l^{-1} 3_l 7_l 1_l^{-1} 17_l
\\
 4_i^{-1} 1_i^{-1} 13_i^{-1} 15_i^{-1} 1_i 4_i 2_i^{-1} 13_i 1_i 3_i^{-1} 17_i^{-1} 1_i 7_i^{-1} 3_i^{-1} 17_i^{-1} 13_i^{-1} 2_i \cdot
 4_j^{-1} 1_j^{-1} 17_j^{-1} 1_j 4_j 2_j^{-1} 13_j 1_j 7_j 4_j 2_j^{-1} 13_j 17_j 3_j 7_j 1_j^{-1} 17_j 3_j 1_j^{-1} 13_j^{-1} 2_j \cdot
 4_k^{-1} 1_k^{-1} 17_k 15_k 13_k 1_k 7_k^{-1} 1_k^{-1} 13_k^{-1} 2_k
\\
 10_i 4_i 2_i^{-1} 13_i 1_i 7_i^{-1} 1_i 7_i^{-1} 1_i^{-1} 13_i^{-1} 2_i 4_i^{-1} 1_i^{-1} 17_i 7_i 1_i^{-1} 13_i^{-1} 2_i 4_i^{-1} 1_i^{-1} 17_i 3_i 10_i^{-1} 4_i^{-1} 15_i \cdot
 10_j^{-1} 15_j^{-1} 4_j 10_j 3_j^{-1} 17_j^{-1} 1_j 4_j 2_j^{-1} 13_j 1_j 7_j^{-1} 17_j^{-1} 1_j 4_j 2_j^{-1} 13_j 1_j 7_j 1_j^{-1} 7_j 1_j^{-1} 13_j^{-1} 2_j 4_j^{-1}
\\
 10_i 1_i^{-1} 13_i^{-1} 2_i 4_i^{-1} 1_i^{-1} 17_i 3_i 10_i^{-1} 4_i^{-1} 15_i 10_i \cdot
 10_j^{-1} 15_j^{-1} 4_j 10_j 3_j^{-1} 17_j^{-1} 1_j 4_j 2_j^{-1} 13_j 1_j 7_j^{-1} 17_j^{-1} 1_j 4_j 2_j^{-1} 13_j 1_j 7_j 1_j^{-1} 7_j 1_j^{-1} 13_j^{-1} 2_j 4_j^{-1} \cdot
 10_k^{-1} 4_k 2_k^{-1} 13_k 1_k 7_k^{-1} 1_k 7_k^{-1} 1_k^{-1} 13_k^{-1} 2_k 4_k^{-1} 1_k^{-1} 17_k 7_k
\\
 4_i 2_i^{-1} 13_i 1_i 7_i 1_i^{-1} 15_i 2_i^{-1} 13_i 17_i 3_i 7_i 1_i 7_i^{-1} 1_i^{-1} 13_i^{-1} 2_i 4_i^{-1} 1_i^{-1} 7_i^{-1} 3_i^{-1} 17_i^{-1} 13_i^{-1} 2_i 15_i^{-1} 1_i \cdot
 4_j 2_j^{-1} 13_j 1_j 7_j 1_j^{-1} 7_j^{-1} 3_j^{-1} 17_j^{-1} 13_j^{-1} 2_j 15_j^{-1} 1_j 7_j^{-1} 1_j^{-1} 13_j^{-1} 2_j 4_j^{-1} 1_j^{-1} 15_j 2_j^{-1} 13_j 17_j 3_j 7_j 1_j
\\
 10_i^{-1} 2_i^{-1} 23_i^{-1} 17_i^{-1} 1_i 7_i^{-1} 3_i^{-1} 17_i^{-1} 1_i 4_i 2_i^{-1} 13_i 1_i 7_i \cdot
 10_j 7_j^{-1} 1_j^{-1} 13_j^{-1} 2_j 4_j^{-1} 1_j^{-1} 17_j 3_j 7_j 1_j^{-1} 17_j 23_j 2_j \cdot
  4_k^{-1} 1_k^{-1} 13_k^{-1} 2_k^{-1} 13_k 17_k 3_k 7_k 1_k 7_k^{-1} 1_k^{-1} 13_k^{-1} 2_k \cdot
  4_l 2_l^{-1} 13_l 1_l 7_l 1_l^{-1} 7_l^{-1} 3_l^{-1} 17_l^{-1} 13_l^{-1} 2_l 13_l 1_l
\\
 4_i 2_i^{-1} 13_i 1_i 7_i 1_i^{-1} 7_i^{-1} 3_i^{-1} 17_i^{-1} 13_i^{-1} 2_i 13_i 1_i \cdot
 4_j^{-1} 1_j^{-1} 13_j^{-1} 15_j^{-1} 1_j 4_j 2_j^{-1} 13_j 1_j 3_j^{-1} 17_j^{-1} 1_j 7_j^{-1} 3_j^{-1} 17_j^{-1} 13_j^{-1} 2_j \cdot
 4_k 2_k^{-1} 13_k 17_k 3_k 7_k 1_k^{-1} 17_k 3_k 1_k^{-1} 13_k^{-1} 2_k 4_k^{-1} 1_k^{-1} 15_k 13_k 1_k \cdot
 4_l^{-1} 1_l^{-1} 13_l^{-1} 2_l^{-1} 13_l 17_l 3_l 7_l 1_l 7_l^{-1} 1_l^{-1} 13_l^{-1} 2_l
\\
 10_i 7_i 1_i^{-1} 17_i 23_i 2_i 13_i 1_i 7_i^{-1} 3_i^{-1} 17_i^{-1} 13_i^{-1} 2_i 4_i^{-1} 1_i^{-1} 13_i^{-1} 2_i^{-1} 13_i 17_i 3_i 7_i 3_i^{-1} 17_i^{-1} 13_i^{-1} 23_i^{-1} 17_i^{-1} 1_i 7_i^{-1} 3_i^{-1} 23_i 2_i \cdot
 10_j^{-1} 2_j^{-1} 23_j^{-1} 3_j 7_j 1_j^{-1} 17_j 23_j 13_j 17_j 3_j 7_j^{-1} 3_j^{-1} 17_j^{-1} 13_j^{-1} 2_j 13_j 1_j 4_j 2_j^{-1} 13_j 17_j 3_j 7_j 1_j^{-1} 13_j^{-1} 2_j^{-1} 23_j^{-1} 17_j^{-1} 1_j 7_j^{-1}
\\
 10_i 7_i 1_i^{-1} 17_i 3_i 10_i^{-1} 2_i^{-1} 13_i 1_i 7_i^{-1} 3_i^{-1} 17_i^{-1} 13_i^{-1} 2_i \cdot
  4_j^{-1} 1_j^{-1} 13_j^{-1} 15_j^{-1} 1_j 4_j 2_j^{-1} 13_j 1_j 3_j^{-1} 17_j^{-1} 1_j 7_j^{-1} 3_j^{-1} 17_j^{-1} 13_j^{-1} 2_j \cdot
 10_k 3_k^{-1} 17_k^{-1} 1_k 7_k^{-1} 10_k^{-1} 2_k^{-1} 13_k 17_k 3_k 7_k 1_k^{-1} 13_k^{-1} 2_k \cdot
  4_l^{-1} 1_l^{-1} 15_l 13_l 1_l 4_l 2_l^{-1} 13_l 17_l 3_l 7_l 1_l^{-1} 17_l 3_l 1_l^{-1} 13_l^{-1} 2_l
\\
 4_i^{-1} 1_i^{-1} 13_i^{-1} 15_i^{-1} 1_i 4_i 2_i^{-1} 13_i 1_i 3_i^{-1} 17_i^{-1} 1_i 7_i^{-1} 3_i^{-1} 17_i^{-1} 13_i^{-1} 2_i \cdot
 4_j^{-1} 1_j^{-1} 15_j 2_j^{-1} 13_j 17_j 3_j 7_j 15_j 13_j 1_j 4_j 2_j^{-1} 13_j 17_j 3_j 7_j 1_j^{-1} 17_j 3_j 1_j^{-1} 13_j^{-1} 2_j 4_j^{-1} 1_j^{-1} 13_j^{-1} 15_j^{-1} 1_j 7_j^{-1} 1_j^{-1} 13_j^{-1} 2_j \cdot
 4_k 2_k^{-1} 13_k 1_k 7_k 1_k^{-1} 15_k 13_k 7_k^{-1} 3_k^{-1} 17_k^{-1} 13_k^{-1} 2_k 15_k^{-1} 1_k
\\
 4_i 2_i^{-1} 13_i 1_i 3_i^{-1} 17_i^{-1} 1_i 7_i^{-1} 3_i^{-1} 17_i^{-1} 13_i^{-1} 2_i 4_i^{-1} 1_i^{-1} 13_i^{-1} 15_i^{-1} 1_i \cdot
 4_j 2_j^{-1} 13_j 1_j 3_j^{-1} 17_j^{-1} 1_j 7_j^{-1} 3_j^{-1} 17_j^{-1} 13_j^{-1} 2_j 4_j^{-1} 1_j^{-1} 13_j^{-1} 15_j^{-1} 17_j 3_j 1_j^{-1} 13_j^{-1} 2_j 4_j^{-1} 1_j^{-1} 15_j 13_j 1_j 4_j 2_j^{-1} 13_j 17_j 3_j 7_j \cdot
 4_k 2_k^{-1} 13_k 17_k 3_k 7_k 1_k^{-1} 17_k 3_k 1_k^{-1} 13_k^{-1} 2_k 4_k^{-1} 1_k^{-1} 15_k 13_k 1_k
\\
 4_i 2_i^{-1} 13_i 1_i 7_i 1_i^{-1} 15_i 2_i^{-1} 13_i 17_i 3_i 7_i 1_i 7_i^{-1} 1_i^{-1} 13_i^{-1} 2_i 4_i^{-1} 1_i^{-1} 7_i^{-1} 3_i^{-1} 17_i^{-1} 13_i^{-1} 2_i 15_i^{-1} 1_i \cdot
 4_j 2_j^{-1} 13_j 1_j 7_j 1_j^{-1} 7_j^{-1} 3_j^{-1} 17_j^{-1} 13_j^{-1} 2_j 15_j^{-1} 17_j^{-1} 1_j 4_j 2_j^{-1} 13_j 1_j 7_j \cdot
 4_k^{-1} 1_k^{-1} 15_k 2_k^{-1} 13_k 17_k 3_k 7_k 1_k 7_k^{-1} 1_k^{-1} 13_k^{-1} 2_k 4_k^{-1} 1_k^{-1} 17_k 1_k 7_k^{-1} 1_k^{-1} 13_k^{-1} 2_k
\\
 4_i^{-1} 1_i^{-1} 17_i 15_i 2_i^{-1} 13_i 17_i 3_i 7_i 1_i 7_i^{-1} 1_i^{-1} 13_i^{-1} 2_i 4_i^{-1} 7_i^{-1} 1_i^{-1} 13_i^{-1} 2_i \cdot
 4_j^{-1} 1_j^{-1} 15_j 2_j^{-1} 13_j 17_j 3_j 7_j 1_j 4_j 2_j^{-1} 13_j 1_j 7_j 1_j^{-1} 7_j^{-1} 3_j^{-1} 17_j^{-1} 13_j^{-1} 2_j 15_j^{-1} 1_j 7_j^{-1} 1_j^{-1} 13_j^{-1} 2_j \cdot
 4_k 2_k^{-1} 13_k 1_k 7_k 1_k^{-1} 7_k^{-1} 3_k^{-1} 17_k^{-1} 13_k^{-1} 2_k 15_k^{-1} 1_k 4_k 2_k^{-1} 13_k 1_k 7_k 1_k^{-1} 17_k^{-1} 1_k
\\
 10_i^{-1} 4_i^{-1} 1_i^{-1} 17_i 3_i 7_i 1_i^{-1} 13_i^{-1} 15_i^{-1} 13_i^{-1} 7_i 10_i^{-1} 15_i^{-1} 4_i 10_i 3_i^{-1} 17_i^{-1} 1_i 4_i 2_i^{-1} 13_i 1_i 4_i 2_i^{-1} 13_i 1_i 7_i^{-1} 17_i^{-1} 1_i 7_i^{-1} 3_i^{-1} 17_i^{-1} 1_i 7_i^{-1} \cdot
 10_j^{-1} 4_j^{-1} 15_j 10_j 7_j^{-1} 13_j 15_j 13_j 1_j 7_j^{-1} 3_j^{-1} 17_j^{-1} 1_j 4_j 10_j 7_j 1_j^{-1} 17_j 3_j 7_j 1_j^{-1} 17_j 7_j 1_j^{-1} 13_j^{-1} 2_j 4_j^{-1} 1_j^{-1} 13_j^{-1} 2_j 4_j^{-1} 1_j^{-1} 17_j 3_j
\\
 4_i^{-1} 1_i^{-1} 13_i^{-1} 15_i^{-1} 1_i 4_i 2_i^{-1} 13_i 1_i 3_i^{-1} 17_i^{-1} 1_i 7_i^{-1} 3_i^{-1} 17_i^{-1} 13_i^{-1} 2_i \cdot
 4_j^{-1} 1_j^{-1} 15_j 2_j^{-1} 13_j 17_j 3_j 7_j 13_j^{-1} 15_j^{-1} 1_j 7_j^{-1} 1_j^{-1} 13_j^{-1} 2_j \cdot
 4_k 2_k^{-1} 13_k 17_k 3_k 7_k 1_k^{-1} 13_k^{-1} 15_k^{-1} 7_k^{-1} 3_k^{-1} 17_k^{-1} 13_k^{-1} 2_k 15_k^{-1} 1_k 4_k 2_k^{-1} 13_k 1_k 7_k 1_k^{-1} 15_k 13_k 15_k 13_k 17_k 3_k 1_k^{-1} 13_k^{-1} 2_k 4_k^{-1} 1_k^{-1} 15_k 13_k 1_k
\\
 10_i^{-1} 15_i^{-1} 4_i 10_i 3_i^{-1} 17_i^{-1} 1_i 4_i 2_i^{-1} 13_i 1_i 7_i^{-1} 15_i^{-1} 4_i 10_i 4_i 2_i^{-1} 13_i 1_i 7_i^{-1} 17_i^{-1} 1_i 7_i^{-1} 3_i^{-1} 17_i^{-1} 1_i 7_i^{-1} 10_i^{-1} 4_i^{-1} 1_i^{-1} 13_i^{-1} 2_i 4_i^{-1} 1_i^{-1} 17_i 3_i 7_i 1_i^{-1} 13_i^{-1} 7_i \cdot
 10_j 7_j 1_j^{-1} 17_j 3_j 7_j 1_j^{-1} 17_j 7_j 1_j^{-1} 13_j^{-1} 2_j 4_j^{-1} 10_j^{-1} 4_j^{-1} 15_j 7_j 1_j^{-1} 13_j^{-1} 2_j 4_j^{-1} 1_j^{-1} 17_j 3_j 10_j^{-1} 4_j^{-1} 15_j 10_j 7_j^{-1} 13_j 1_j 7_j^{-1} 3_j^{-1} 17_j^{-1} 1_j 4_j 2_j^{-1} 13_j 1_j 4_j
\\
 10_i^{-1} 2_i^{-1} 13_i 17_i 3_i 7_i 1_i^{-1} 13_i^{-1} 23_i^{-1} 7_i 1_i^{-1} 13_i^{-1} 2_i 4_i^{-1} 10_i^{-1} 15_i^{-1} 4_i 10_i 3_i^{-1} 17_i^{-1} 1_i 4_i 2_i^{-1} 13_i 1_i 4_i 2_i^{-1} 13_i 1_i 7_i^{-1} 17_i^{-1} 1_i 7_i^{-1} 3_i^{-1} 17_i^{-1} 1_i 7_i^{-1} \cdot
 10_j 7_j 1_j^{-1} 17_j 3_j 7_j 1_j^{-1} 17_j 7_j 1_j^{-1} 13_j^{-1} 2_j 4_j^{-1} 10_j^{-1} 4_j^{-1} 15_j 7_j 4_j 2_j^{-1} 13_j 1_j 7_j^{-1} 23_j 13_j 1_j 7_j^{-1} 3_j^{-1} 17_j^{-1} 13_j^{-1} 2_j \cdot
 10_k 1_k^{-1} 13_k^{-1} 2_k 4_k^{-1} 1_k^{-1} 17_k 3_k 10_k^{-1} 4_k^{-1} 15_k 10_k 7_k^{-1} 15_k^{-1} 4_k
\\
 10_i^{-1} 15_i^{-1} 4_i 10_i 3_i^{-1} 17_i^{-1} 1_i 4_i 2_i^{-1} 13_i 1_i 4_i 2_i^{-1} 13_i 1_i 7_i^{-1} 17_i^{-1} 1_i 7_i^{-1} 3_i^{-1} 17_i^{-1} 1_i 7_i^{-1} 10_i^{-1} 4_i^{-1} 1_i^{-1} 17_i 3_i 7_i 1_i^{-1} 17_i 3_i 1_i^{-1} 17_i 3_i 7_i 1_i^{-1} 13_i^{-1} 7_i \cdot
 10_j 7_j^{-1} 15_j^{-1} 4_j 10_j 1_j^{-1} 13_j^{-1} 2_j 4_j^{-1} 1_j^{-1} 17_j 3_j 10_j^{-1} 4_j^{-1} 15_j \cdot
 10_k 7_k 1_k^{-1} 17_k 3_k 7_k 1_k^{-1} 17_k 7_k 1_k^{-1} 13_k^{-1} 2_k 4_k^{-1} 10_k^{-1} 4_k^{-1} 15_k 13_k 1_k 7_k^{-1} 3_k^{-1} 17_k^{-1} 1_k 3_k^{-1} 17_k^{-1} 1_k 7_k^{-1} 3_k^{-1} 17_k^{-1} 1_k 4_k
\\
 4_i^{-1} 1_i^{-1} 17_i 15_i 13_i 1_i 7_i^{-1} 1_i^{-1} 13_i^{-1} 2_i \cdot
 4_j 2_j^{-1} 13_j 17_j 3_j 7_j 1_j^{-1} 17_j 3_j 1_j^{-1} 13_j^{-1} 2_j 4_j^{-1} 1_j^{-1} 15_j 13_j 1_j \cdot
 4_k^{-1} 1_k^{-1} 15_k 13_k 1_k 4_k 2_k^{-1} 13_k 17_k 3_k 7_k 1_k^{-1} 17_k 3_k 1_k^{-1} 13_k^{-1} 2_k \cdot
 4_l 2_l^{-1} 13_l 1_l 3_l^{-1} 17_l^{-1} 13_l^{-1} 2_l 4_l^{-1} 1_l^{-1} 13_l^{-1} 15_l^{-1} 1_l 4_l 2_l^{-1} 13_l 1_l 3_l^{-1} 17_l^{-1} 13_l^{-1} 15_l^{-1} 1_l 7_l^{-1} 3_l^{-1} 17_l^{-1} 13_l^{-1} 2_l 4_l^{-1} 1_l^{-1} 13_l^{-1} 15_l^{-1} 1_l 4_l 2_l^{-1} 13_l 1_l 3_l^{-1} 17_l^{-1} 17_l^{-1} 1_l
\\
 10_i 7_i^{-1} 13_i 1_i 7_i^{-1} 3_i^{-1} 17_i^{-1} 1_i 4_i 10_i 7_i 1_i^{-1} 17_i 3_i 7_i 1_i^{-1} 17_i 23_i 13_i 1_i 7_i^{-1} 3_i^{-1} 17_i^{-1} 13_i^{-1} 2_i 4_i^{-1} 1_i^{-1} 17_i 3_i 7_i 1_i^{-1} 13_i^{-1} 7_i 1_i^{-1} 13_i^{-1} 2_i 4_i^{-1} 1_i^{-1} 17_i 3_i 10_i^{-1} 4_i^{-1} 15_i \cdot
 10_j^{-1} 4_j^{-1} 1_j^{-1} 17_j 3_j 7_j 1_j^{-1} 13_j^{-1} 7_j 10_j^{-1} 15_j^{-1} 4_j 10_j 3_j^{-1} 17_j^{-1} 1_j 4_j 2_j^{-1} 13_j 1_j 7_j^{-1} 13_j 1_j 7_j^{-1} 3_j^{-1} 17_j^{-1} 1_j 4_j 2_j^{-1} 13_j 17_j 3_j 7_j 1_j^{-1} 13_j^{-1} 23_j^{-1} 17_j^{-1} 1_j 7_j^{-1} 3_j^{-1} 17_j^{-1} 1_j 7_j^{-1}
\\
 4_i^{-1} 1_i^{-1} 15_i 2_i^{-1} 13_i 17_i 3_i 7_i 7_i 1_i^{-1} 13_i^{-1} 2_i 4_i^{-1} 7_i^{-1} 7_i^{-1} 1_i 4_i 2_i^{-1} 13_i 17_i 3_i 7_i 1_i^{-1} 17_i 3_i 1_i^{-1} 13_i^{-1} 2_i \cdot
 4_j 2_j^{-1} 13_j 1_j 3_j^{-1} 17_j^{-1} 1_j 7_j^{-1} 3_j^{-1} 17_j^{-1} 13_j^{-1} 2_j 4_j^{-1} 1_j^{-1} 13_j^{-1} 15_j^{-1} 1_j \cdot
 4_k 2_k^{-1} 13_k 1_k 7_k^{-1} 7_k^{-1} 3_k^{-1} 17_k^{-1} 13_k^{-1} 2_k 15_k^{-1} 1_k 4_k 2_k^{-1} 13_k 1_k 3_k^{-1} 17_k^{-1} 1_k 7_k^{-1} 3_k^{-1} 17_k^{-1} 13_k^{-1} 2_k 4_k^{-1} 1_k^{-1} 7_k 7_k \cdot
 4_l^{-1} 1_l^{-1} 15_l 13_l 1_l 4_l 2_l^{-1} 13_l 17_l 3_l 7_l 1_l^{-1} 17_l 3_l 1_l^{-1} 13_l^{-1} 2_l
\\
 10_i^{-1} 15_i^{-1} 4_i 10_i 3_i^{-1} 17_i^{-1} 1_i 4_i 2_i^{-1} 13_i 1_i 4_i 2_i^{-1} 13_i 1_i 7_i^{-1} 17_i^{-1} 1_i 7_i^{-1} 3_i^{-1} 17_i^{-1} 1_i 7_i^{-1} 10_i^{-1} 2_i^{-1} 13_i 17_i 3_i 7_i 1_i^{-1} 13_i^{-1} 23_i^{-1} 17_i^{-1} 1_i 7_i^{-1} 3_i^{-1} 17_i^{-1} 1_i 7_i^{-1} 10_i^{-1} 2_i^{-1} 13_i 1_i 4_i 10_i 7_i 1_i^{-1} 17_i 3_i 7_i 1_i^{-1} 17_i 7_i 1_i^{-1} 13_i^{-1} 2_i 4_i^{-1} 1_i^{-1} 13_i^{-1} 2_i 4_i^{-1} 1_i^{-1} 17_i 3_i 10_i^{-1} 4_i^{-1} 15_i 10_i 7_i^{-1} 13_i 1_i 7_i^{-1} 3_i^{-1} 17_i^{-1} 13_i^{-1} 2_i 10_i 7_i 1_i^{-1} 17_i 3_i 7_i 1_i^{-1} 17_i 23_i 13_i 1_i 7_i^{-1} 3_i^{-1} 17_i^{-1} 13_i^{-1} 2_i 4_i^{-1} 1_i^{-1} 17_i 3_i 7_i 1_i^{-1} 13_i^{-1} 7_i
\\
 10_i^{-1} 15_i^{-1} 4_i 10_i 3_i^{-1} 17_i^{-1} 1_i 4_i 2_i^{-1} 13_i 1_i 4_i 2_i^{-1} 13_i 1_i 7_i^{-1} 17_i^{-1} 1_i 7_i^{-1} 3_i^{-1} 17_i^{-1} 1_i 7_i^{-1} 10_i^{-1} 4_i^{-1} 1_i^{-1} 13_i^{-1} 2_i 4_i^{-1} 1_i^{-1} 17_i 3_i 7_i 1_i^{-1} 13_i^{-1} 7_i 4_i 2_i^{-1} 13_i 1_i 7_i^{-1} 1_i 7_i^{-1} 1_i^{-1} 13_i^{-1} 2_i 4_i^{-1} 1_i^{-1} 17_i 7_i \cdot
 10_j^{-1} 4_j^{-1} 15_j 10_j 7_j^{-1} 13_j 1_j 7_j^{-1} 3_j^{-1} 17_j^{-1} 1_j 4_j 2_j^{-1} 13_j 1_j 4_j 10_j 7_j 1_j^{-1} 17_j 3_j 7_j 1_j^{-1} 17_j 7_j 1_j^{-1} 13_j^{-1} 2_j 4_j^{-1} 7_j^{-1} 17_j^{-1} 1_j 4_j 2_j^{-1} 13_j 1_j 7_j 1_j^{-1} 7_j 1_j^{-1} 13_j^{-1} 2_j 4_j^{-1} 1_j^{-1} 13_j^{-1} 2_j 4_j^{-1} 1_j^{-1} 17_j 3_j
\\
 10_i 7_i 1_i^{-1} 17_i 3_i 7_i 1_i^{-1} 17_i 7_i 1_i^{-1} 13_i^{-1} 2_i 4_i^{-1} 7_i 10_i^{-1} 2_i^{-1} 23_i^{-1} 17_i^{-1} 1_i 7_i^{-1} 10_i^{-1} 4_i^{-1} 15_i 10_i 7_i^{-1} 13_i 1_i 7_i^{-1} 3_i^{-1} 17_i^{-1} 1_i 4_i 2_i^{-1} 13_i 17_i 3_i 7_i 1_i^{-1} 13_i^{-1} 23_i^{-1} 17_i^{-1} 1_i 7_i^{-1} 3_i^{-1} 17_i^{-1} 1_i 7_i^{-1} 10_i^{-1} 2_i^{-1} 13_i 1_i 4_i \cdot
 10_j^{-1} 15_j^{-1} 4_j 10_j 7_j 1_j^{-1} 17_j 23_j 2_j 10_j 7_j^{-1} 4_j 2_j^{-1} 13_j 1_j 7_j^{-1} 17_j^{-1} 1_j 7_j^{-1} 3_j^{-1} 17_j^{-1} 1_j 7_j^{-1} 10_j^{-1} 4_j^{-1} 1_j^{-1} 13_j^{-1} 2_j 10_j 7_j 1_j^{-1} 17_j 3_j 7_j 1_j^{-1} 17_j 23_j 13_j 1_j 7_j^{-1} 3_j^{-1} 17_j^{-1} 13_j^{-1} 2_j 4_j^{-1} 1_j^{-1} 17_j 3_j 7_j 1_j^{-1} 13_j^{-1} 7_j
\\
 10_i^{-1} 4_i^{-1} 15_i 10_i 4_i 2_i^{-1} 13_i 1_i 7_i^{-1} 1_i 7_i^{-1} 1_i^{-1} 13_i^{-1} 2_i 4_i^{-1} 1_i^{-1} 17_i 7_i 1_i^{-1} 13_i^{-1} 2_i 4_i^{-1} 1_i^{-1} 17_i 3_i \cdot
 10_j^{-1} 15_j^{-1} 1_j 4_j 2_j^{-1} 13_j 1_j 3_j^{-1} 17_j^{-1} 13_j^{-1} 7_j 10_j^{-1} 15_j^{-1} 4_j 10_j 3_j^{-1} 17_j^{-1} 1_j 4_j 2_j^{-1} 13_j 1_j 4_j 2_j^{-1} 13_j 1_j 7_j^{-1} 17_j^{-1} 1_j 7_j^{-1} 3_j^{-1} 17_j^{-1} 1_j 7_j^{-1} 10_j^{-1} 4_j^{-1} 1_j^{-1} 17_j 3_j 7_j 1_j^{-1} 13_j^{-1} 7_j \cdot
 10_k 7_k^{-1} 13_k 1_k 7_k^{-1} 3_k^{-1} 17_k^{-1} 1_k 4_k 10_k 7_k 1_k^{-1} 17_k 3_k 7_k 1_k^{-1} 17_k 7_k 1_k^{-1} 13_k^{-1} 2_k 4_k^{-1} 7_k^{-1} 17_k^{-1} 1_k 4_k 2_k^{-1} 13_k 1_k 7_k 1_k^{-1} 7_k 1_k^{-1} 13_k^{-1} 2_k 4_k^{-1} 7_k^{-1} 13_k 17_k 3_k 1_k^{-1} 13_k^{-1} 2_k 4_k^{-1} 1_k^{-1} 15_k
\\
 10_i^{-1} 15_i^{-1} 4_i 10_i 3_i^{-1} 17_i^{-1} 1_i 4_i 2_i^{-1} 13_i 1_i 4_i 2_i^{-1} 13_i 1_i 7_i^{-1} 17_i^{-1} 1_i 7_i^{-1} 3_i^{-1} 17_i^{-1} 1_i 7_i^{-1} 10_i^{-1} 2_i^{-1} 13_i 17_i 3_i 7_i 1_i^{-1} 13_i^{-1} 23_i^{-1} 1_i 7_i^{-1} 1_i^{-1} 13_i^{-1} 2_i 4_i^{-1} 1_i^{-1} 17_i 7_i 1_i^{-1} 13_i^{-1} 2_i 4_i^{-1} 7_i^{-1} 13_i 1_i 7_i^{-1} 3_i^{-1} 17_i^{-1} 13_i^{-1} 2_i 4_i^{-1} 1_i^{-1} 15_i 13_i 1_i 4_i 2_i^{-1} 13_i 17_i 3_i 7_i 1_i^{-1} 13_i^{-1} 7_i \cdot
 10_j 7_j 1_j^{-1} 17_j 3_j 7_j 1_j^{-1} 17_j 7_j 1_j^{-1} 13_j^{-1} 2_j 4_j^{-1} 1_j^{-1} 13_j^{-1} 2_j 4_j^{-1} 1_j^{-1} 17_j 3_j 10_j^{-1} 4_j^{-1} 15_j 10_j 4_j 2_j^{-1} 13_j 1_j 7_j^{-1} 17_j^{-1} 1_j 4_j 2_j^{-1} 13_j 1_j 7_j 1_j^{-1} 23_j 13_j 1_j 7_j^{-1} 3_j^{-1} 17_j^{-1} 13_j^{-1} 2_j 4_j^{-1} 1_j^{-1} 13_j^{-1} 2_j 4_j^{-1} 1_j^{-1} 13_j^{-1} 15_j^{-1} 1_j 4_j 2_j^{-1} 13_j 1_j 4_j
\\
 10_i^{-1} 2_i^{-1} 13_i 1_i 4_i 2_i^{-1} 13_i 17_i 3_i 7_i 1_i^{-1} 13_i^{-1} 23_i^{-1} 17_i^{-1} 1_i 7_i^{-1} 3_i^{-1} 17_i^{-1} 1_i 7_i^{-1} 10_i^{-1} 2_i^{-1} 13_i 17_i 3_i 7_i 1_i^{-1} 13_i^{-1} 23_i^{-1} 1_i 7_i^{-1} 1_i^{-1} 13_i^{-1} 2_i 4_i^{-1} 1_i^{-1} 17_i 13_i 1_i 7_i^{-1} 3_i^{-1} 17_i^{-1} 1_i 4_i 2_i^{-1} 13_i 17_i 3_i 7_i 1_i^{-1} 13_i^{-1} 23_i^{-1} 17_i^{-1} 1_i 7_i^{-1} 3_i^{-1} 17_i^{-1} 1_i 7_i^{-1} \cdot
 10_j 7_j 1_j^{-1} 17_j 3_j 7_j 1_j^{-1} 17_j 23_j 13_j 1_j 7_j^{-1} 3_j^{-1} 17_j^{-1} 13_j^{-1} 2_j 4_j^{-1} 1_j^{-1} 13_j^{-1} 2_j 4_j^{-1} 1_j^{-1} 17_j 3_j 7_j 1_j^{-1} 13_j^{-1} 7_j 4_j 2_j^{-1} 13_j 1_j 7_j^{-1} 23_j 13_j 1_j 7_j^{-1} 3_j^{-1} 17_j^{-1} 13_j^{-1} 2_j \cdot
 10_k 7_k 1_k^{-1} 17_k 3_k 7_k 1_k^{-1} 17_k 23_k 13_k 1_k 7_k^{-1} 3_k^{-1} 17_k^{-1} 13_k^{-1} 2_k 4_k^{-1} 1_k^{-1} 17_k 3_k 7_k 1_k^{-1} 13_k^{-1} 17_k^{-1} 1_k 4_k 2_k^{-1} 13_k 1_k 7_k 1_k^{-1} 7_k 1_k^{-1} 13_k^{-1} 2_k 4_k^{-1} 7_k^{-1} 13_k 1_k 7_k^{-1} 3_k^{-1} 17_k^{-1} 1_k 4_k
\end{eqnarray*}
\forgotten

\subsection{Breaking $N$ into sections}\label{modcent}

The subgroup $N$ defined above is generated, as a normal subgroup
of $F$ invariant under the action of $S_n$, by $49$ elements. Let
$K = F/N$, so that $$\pitil^{\Aff} \isom S_n \semidirect K.$$ %
For convenience, we would like to work in $F^{*}$ rather than $F$
(see Subsection \ref{ss:struc} for the definitions). It is easy to
see that a normal subgroup of $F$ which is generated by elements
with at least one trivial entry in $F^* = F_1 \times \dots \times
F_n$, is
also normal in $F^*$. In particular $N$ is normal in $F^*$. %
We now define 
\begin{equation}\label{Kstardef}
K^*= F^*/N. \end{equation} %
Since $N \leq F$, the map $\ab \co K^* \ra \Z^{10}$ is well
defined, and has $K$ as its kernel. As in $F$, $K$ is generated by
the elements $\omega_\alpha\omega_\beta^{-1}$ (i.e.\ their images
in $K^*$).

The images of a word in $F^* = F_1 \times\dots\times F_{18}$ under
the projections to the $F_i$ are called `sections'; e.g. the
sections of $u_i w_j v_k$ are $u_i,w_j$ and $v_k$ (where $u_i \in
F_i$, $w_j \in F_j$ and $v_k \in F_k$). Let $\hat{N}$ be the group
generated by sections of elements in $N$. Alternatively, $\hat{N}$
is generated (as a normal and invariant subgroup of $F^*$) by the
$122 = 1 \cdot 1 + 29 \cdot 2 + 13 \cdot 3 + 6 \cdot 4$ sections
of the $49$ generators of $N$. Obviously $N \sub \hat{N}$.

\begin{prop}\label{Nhat}
$\hat{N}/N$ is central in $F^*/N$. Moreover, $\hat{N}/N$ is fixed
(element-wise) under the action of $S_n$.
\end{prop}
\begin{proof}
The first claim holds for any normal subgroup of $F^*$. Indeed,
let $u = u_1 u_2 \dots u_{n}$ be an element of $N$ (with $u_i \in
F_i$), and consider the section $u_i$. By definition of $F^*$,
$u_i$ commutes with every generator $\omega_j$ for $j \neq i$. But
at the same time, $u_i$ can be expressed (modulo $N$) as a product
of the other sections of $u$, and so it also commutes with the
generators $\omega_i$.

For the second claim we need to know that $N$ is invariant and
generated by words with less than $n$ non-trivial sections (this
is indeed the case in $N$, where the generators have at most $4$
sections). Without loss of generality let $u_1 u_2 \dots u_k \in
N$ where $k < n$, then $\tau(u_1) u_2 \dots u_k \in N$ for $\tau$
the transposition $(1\, n)$, so modulo $N$ we have that $\tau(u_1)
\equiv u_1$. Acting now with an arbitrary permutation $\s \in S_n$
which fixes the index $1$, we get that $\s(u_1) \equiv u_1$ and so
$u_1$ is fixed under $S_n$ (as an element of $\hat{N}/N$).
\end{proof}

By definition, every section of an element of $\hat{N}$ is also in
$\hat{N}$, so that $\hat{N} = (F_1 \cap \hat{N})\cdots (F_n \cap
\hat{N})$. It follows that %
\begin{equation}\label{H1Hn}
F^* / \hat{N} = (F_1
\times \dots \times F_n) / \hat{N} \isom H_1 \times \dots \times
H_n,
\end{equation} where $H_i = F_i / (F_i \cap \hat{N})$.
Moreover, the groups $H_i$ are all isomorphic, being permuted by
$S_n$.

\ifXY
\begin{figure}
\begin{equation*}
\xymatrix@C=12pt@R=10pt{ 
          {}
        & {}
        & {}
        & {} F^* \ar@{-}[ld]
\\
          {}
        & {}
        & {} F_i \hat{N} \ar@{-}[ld] \ar@{-}[rd]^{H_i}
        & {}
\\
          {}
        & {} F_i N \ar@{-}[ld] \ar@{-}[rd]^{H_i}
        & {}
        & {} \hat{N} \ar@{-}[ld]
\\
          {} F_i \ar@{-}[rd]^{H_i}
        & {}
        & {} N \cdot F_i \cap \hat{N} \ar@{-}[ld] \ar@{-}[rd]
        & {}
\\
          {}
        & {} F_i \cap \hat{N} \ar@{-}[rd]
        & {}
        & {} N \ar@{-}[ld]
\\
          {}
        & {}
        & {} F_i \cap N \ar@{-}[rd]
        & {}
\\
    {}
        & {}
        & {}
        & {} 1
}
\end{equation*}
\caption{Decomposition of $F^*/\hat{N}$}
\end{figure}
\else
\fi 

As an illustration, recall the generator of $N$ from \Eq{rel}. The
first section is
$$1_i^{-1} \sp 13_i^{-1} \sp 2_i \sp 4_i^{-1} \sp 7_i \sp 4_i \sp 2_i^{-1} \sp 13_i
\sp 1_i \sp 7_i^{-1} \in \hat{N},$$ which can written shortly as
$$1^{-1} \cdot 13^{-1} \cdot 2 \cdot 4^{-1}\cdot 7\cdot 4\cdot 2^{-1}\cdot 13\cdot
1\cdot 7^{-1},$$ to represent an element of $(N \cdot F_i \cap
\hat{N})/N$ for arbitrary $i=1,\dots,n$.

\subsection{The group $H$}\label{groupH}

The decomposition of \Eq{H1Hn} presents $K^* = F^*/N$ modulo its
central subgroup $\hat{N}/N$ as a product of $n = 18$ isomorphic
groups, which from now on we denote by $H$.  This group is
generated by $\Omega = \set{1,2,\dots,17,23}$ with $122$ relations
obtained from the sections of the previous $49$ generators of $N$.
It turns out that the $122$ relations include some repetition, so
we only have $97$ relations. Here are the shortest seven of them:
\def\spp{{\,\cdot\,}}
\begin{eqnarray*}
10\spp 7\spp 10^{-1}\spp 4^{-1}\spp 15 & = & e \\
17^{-1}\spp 1\spp 3^{-1}\spp 17^{-1}\spp 13^{-1}\spp 15^{-1}\spp 1\spp 7^{-1}\spp 3^{-1} & = & e \\
17\spp 3\spp 1^{-1}\spp 17\spp 3\spp 7\spp 1^{-1}\spp 15\spp 13 & = & e \\
4^{-1}\spp 1^{-1}\spp 17\spp 15\spp 13\spp 1\spp 7^{-1}\spp 1^{-1}\spp 13^{-1}\spp 2 & = & e \\
4^{-1}\spp 7^{-1}\spp 4\spp 2^{-1}\spp 13\spp 1\spp 7\spp 1^{-1}\spp 13^{-1}\spp 2 & = & e \\
4^{-1}\spp 7\spp 4\spp 2^{-1}\spp 13\spp 1\spp 7^{-1}\spp 1^{-1}\spp 13^{-1}\spp 2 & = & e \\
4\spp 2^{-1}\spp 13\spp 1\spp 7\spp 1^{-1}\spp 13^{-1}\spp 15^{-1}
\spp 17^{-1}\spp 1 & = & e
\end{eqnarray*}

At this point we apply the shortening algorithm described in
Subsection \ref{ss:short}, and obtain $67$ relations on the ten
generators, of total length $382$. 
See [App4] (see the Appendix) for the complete list. Some of the
relations are given below:
\begin{eqnarray*}
  10 \spp 17^{-1}\spp 7  & = & e\\
  10 \spp 7 \spp 17^{-1} & = & e\\
  13 \spp 23^{-1} \spp 15 & = & e\\
  4 \spp 7^{-1} \spp 15^{-1} & = & e\\
  4\spp15^{-1}\spp7^{-1}  & = & e\\
\forget
  7^{-1} 10^{-1} 17 & = & e\\
  10 3 10^{-1} 3^{-1} & = & e\\
  13 7 13^{-1} 7^{-1} & = & e\\
  13^{-1} 4^{-1} 3^{-1} 2 & = & e\\
  2 3^{-1} 7^{-1} 23^{-1} & = & e\\
  23 1 3^{-1} 10^{-1} & = & e\\
  23 10^{-1} 23^{-1} 10 & = & e\\
  23 17 23^{-1} 17^{-1} & = & e\\
  23 17^{-1} 23^{-1} 17 & = & e\\
  23^{-1} 17^{-1} 23 17 & = & e\\
  23^{-1} 3 1^{-1} 10 & = & e\\
  23^{-1} 7^{-1} 23 7 & = & e\\
  2^{-1} 3 7 23 & = & e\\
  4 13 3 2^{-1} & = & e\\
  7 13^{-1} 7^{-1} 13 & = & e\\
  7 23 2^{-1} 3 & = & e\\
  7^{-1} 23^{-1} 2 3^{-1} & = & e\\
  13 7^{-1} 3^{-1} 15 2 & = & e\\
  15^{-1} 1 13^{-1} 1^{-1} 23^{-1} & = & e\\
  17^{-1} 23^{-1} 1 7^{-1} 2^{-1} & = & e\\
  23^{-1} 3^{-1} 15^{-1} 3 13^{-1} & = & e\\
\forgotten
  7\spp23\spp2\spp17\spp1^{-1} & = & e\\
  7^{-1}\spp23^{-1}\spp2^{-1}\spp1\spp17^{-1} & = & e\\
  1\spp7\spp10^{-1}\spp4\spp2^{-1}\spp4^{-1} & = & e.\\
\forget
  10^{-1} 3 13 4 10 2^{-1} & = & e\\
  17 23 3^{-1} 17^{-1} 10^{-1} 1 & = & e\\
  2^{-1} 7 10^{-1} 13 1 13^{-1} & = & e\\
  3^{-1} 23 2 4 7 13 & = & e\\
  4 2 4^{-1} 13^{-1} 2^{-1} 13 & = & e\\
  7^{-1} 13^{-1} 10 2 13 1^{-1} & = & e\\
  7^{-1} 1^{-1} 13^{-1} 2 13 10 & = & e\\
  7^{-1} 2^{-1} 13^{-1} 3 4^{-1} 23 & = & e\\
  17^{-1} 23^{-1} 13^{-1} 3 10 2^{-1} 4^{-1} & = & e\\
  4^{-1} 1 7^{-1} 1^{-1} 17 4 10^{-1} & = & e\\
  10^{-1} 13 10 15 2^{-1} 15^{-1} 13^{-1} 2 & = & e\\
  10^{-1} 23 2^{-1} 10^{-1} 1 13 17 15 & = & e\\
  1^{-1} 23^{-1} 1 7 1^{-1} 7^{-1} 1 23 & = & e\\
  2 7^{-1} 2^{-1} 13 17 15 23^{-1} 10^{-1} & = & e\\
  23 10 2 4^{-1} 1^{-1} 23^{-1} 13^{-1} 23 & = & e\\
  23 17 2 13 15 23^{-1} 7 1^{-1} & = & e\\
  23 7 2 10^{-1} 2^{-1} 10 4^{-1} 13^{-1} & = & e\\
  23 7^{-1} 10 15^{-1} 3^{-1} 15 10^{-1} 2 & = & e\\
  2^{-1} 3 13^{-1} 10 4^{-1} 13^{-1} 3 23 4^{-1} 1^{-1} & = & e\\
  7^{-1} 3^{-1} 17^{-1} 13^{-1} 17 3 17^{-1} 15^{-1} 10 7^{-1} & = & e\\
  13 7^{-1} 2^{-1} 15 2^{-1} 10^{-1} 1 4^{-1} 3 4 23^{-1} & = & e\\
  17 23 2 10 4 2^{-1} 23^{-1} 7 10^{-1} 10^{-1} 13 & = & e\\
  13^{-1} 4^{-1} 4^{-1} 7^{-1} 2 7 10^{-1} 7 1^{-1} 13^{-1} 10 10 & = & e\\
  13^{-1} 4^{-1} 17^{-1} 4^{-1} 23^{-1} 10^{-1} 10^{-1} 2^{-1} 13^{-1} 10 10 7^{-1} 3 10 7^{-1} & = & e\\
  17 23 10 7^{-1} 4^{-1} 1 10^{-1} 4^{-1} 13^{-1} 10^{-1} 7 7 1^{-1} 13^{-1} 23 & = & e\\
  23 17 2 1^{-1} 13^{-1} 2^{-1} 4^{-1} 13^{-1} 3 15^{-1} 13^{-1} 17 7 10^{-1} 17^{-1} 7^{-1} 10 13 & = & e
\forgotten
\end{eqnarray*}

Substituting some of the generators in terms of others and
simplifying the relations further, we obtain the following
equivalent description of $H$. Set $c = 13\cdot 4$, so the
generator $13$ can be replaced by $c \cdot 4^{-1}$. Then we obtain
the following equalities in $H$:
\begin{eqnarray}
15 & = & 7^{-1} \cdot 4, \nonumber \\
17 & = & 10 \cdot 7, \nonumber \\
23 & = & c \cdot 7^{-1}, \label{equalH} \\ 
2 & = & c^2 \cdot 7^{-1} \cdot 10^{-1} \cdot 1, \nonumber\\ 
3 & = & c \cdot 7^{-1} \cdot 10^{-1} \cdot  1. \nonumber
\end{eqnarray}
In particular $H = \sg{1,4,7,10,c}$. The shortening algorithm
yields the following presentation on these five generators:
\begin{equation}\label{defrelsH}
\begin{matrix}
{}[1,c] & = & [7,c] & = & [4,c] & = & [10,c] & = & e,\\
{}      &  & [10,1] & = & [10,7] & = & [4,7] & = & e,\\
{}      &  &        &   & [1,7]^{-1} & = & [4,10] & = & c^3,\\
{}      &  &        &   &       &   &  [4,1] & = & c^{-2} 7^{2}.
\end{matrix}
\end{equation}
\begin{cor}\label{Hisnilp}
$H$ is nilpotent of class $3$.
\end{cor}
\begin{proof}
The commutator subgroup of $H$ is generated by $c^3$ and $c^{-2}
7^2$, so $H' = \sg{c,7^2}$ (isomorphic to $\Z^2$). It follows that
$[H,H'] = [H,\sg{7^2}^H] = [\sg{1}^H,\sg{7^2}^H] = \sg{c^6}$ which
is central.
\end{proof}

We note that the subgroup $\sg{c,7,10}$ is isomorphic to $\Z^3$,
and that $H/\sg{c,7,10} \isom \Z^2$. Thus $H$ is an extension
$$1 \llra \Z^3 \llra H \llra \Z^2 \llra 1,$$
and so $H$ has derived length $2$.

\forget%
Let $H_1 = H/\Cent(H)$ and $H_2 = H_1 / \Cent(H_1)$. It is
easy to that $\Cent(H) = \sg{c} \isom \Z$. Therefore, $H_1$ is
generated by $1,4,7,10$ subject to the relations
\begin{eqnarray*}
7,10 & \in & \Cent(H_1),\\
{}[1,4] & = & 7^{-2}.
\end{eqnarray*}
Next, $\Cent(H_1) = \sg{7,10} \isom \Z^2$, and so $H_2 = \sg{1, 4
\subjectto [1,4] = e} \isom \Z^2$. %
\forgotten

\subsection{Simplifying $N$}

In Subsection \ref{modcent} we proved that $F^*/\hat{N}$
decomposes in a natural way into a product $H_1 \times \dots
\times H_{18}$, with each $H_i$ isomorphic to the group $H$
studied in Subsection \ref{groupH}. This group is generated by the
$10$ generators $\Omega$, with $15$ defining relations (out of
which, five relations express redundant generators, and the other
ten are related to commutators of the remaining generators).

We lift this description back to $K^* = F^* / N$, recalling that
$N$ is generated by $49$ relations (as a normal subgroup of $F^*$
which is invariant under $S_{18}$). In order to simplify these
relations, we apply the defining relations of $H$ --- this time
not as relations, but by equating each one of them to a newly
defined element of $\hat{N}/N \sub \Cent(K^*)$. In other words
(after replacing $13_i$ by $c_i 4^{-1}_i$ throughout), we define
for $i = 1,\dots,18$, $\theta_i^{(15)} = 15_i 4_i^{-1} 7_i$,
$\theta_i^{(17)} = 17_i 7_i^{-1} 10_i^{-1}$, and likewise
$\theta_i^{(23)}$, $\theta_i^{(2)}$ and $\theta_i^{(3)}$ (see
\Eq{equalH}). In a similar manner we set $\theta_i^{(1,c)} =
[1_i,c_i]$, $\theta_i^{(7,c)} = [7_i,c_i]$ etc., up to
$\theta_i^{(4,1)} = [4_i,1_i]7_i^{-2}c_i^2$ (see \Eq{defrelsH}).
The $18\cdot 10$ equations of the form $[1_i,c_i] =
\theta_i^{(1,c)}$ are called the commutator relations of $K^*$.
Since modulo $\hat{N}$ all these products become trivial, all the
$18\cdot 15$ new elements $\theta_i^{(\cdot)}$ just defined,
belong to $\hat{N}/N$. It follows from Proposition \ref{Nhat} that
the $\theta_i^{(\cdot)}$ are in fact independent of the index $i$,
being fixed under $S_n$. Nevertheless, since we anyway have to
simplify the $49$ relations, we keep the indices to help tracking
the computation.

The first step in simplifying the generators of $N$ is to express
the `redundant generators' $15_i,17_i,23_i,2_i$ and $3_i$ in terms
of $1_i,4_i,7_i,10_i,c_i$ and the variables $\theta_i^{(\cdot)}$.
For example, the relation given in \Eq{rel} becomes
\begin{eqnarray*}
& & 1_i^{-1} 4_i c_i^{-1} \theta_i^{(2)} c_i^2 7_i^{-1} 10_i^{-1}
1_i 4_i^{-1} 7_i 4_i 1_i^{-1} 10_i 7_i c_i^{-2}
{\theta_i^{(2)}}^{-1}
c_i 4^{-1}_i 1_i 7_i^{-1} \\
& & {} \quad\cdot \, 7_j 1_j^{-1} 4_j c_j^{-1} \theta_j^{(2)}
c_j^2 7_j^{-1} 10_j^{-1} 1_j 4_j^{-1} 7_j^{-1} 4_j 1_j^{-1} 10_j
7_j c_j^{-2} {\theta_j^{(2)}}^{-1} c_j 4^{-1}_j
1_j = e, %
\end{eqnarray*}
and then, since $\theta_i^{(2)}$ is central,
\begin{eqnarray}\label{rel+}
\begin{array}{cc}
 & 1_i^{-1} 4_i c_i 7_i^{-1} 10_i^{-1} 1_i 4_i^{-1} 7_i 4_i
1_i^{-1} 10_i 7_i c_i^{-1} 4^{-1}_i 1_i 7_i^{-1} \\
 & {} \qquad\qquad\cdot \, 7_j 1_j^{-1} 4_j c_j 7_j^{-1} 10_j^{-1} 1_j
4_j^{-1} 7_j^{-1} 4_j 1_j^{-1} 10_j 7_j c_j^{-1} 4^{-1}_j
1_j = e %
\end{array}
\end{eqnarray}
(it so happens that almost all the relations are balanced in terms
the $\theta_i^{(\cdot)}$, and so they rarely show up in the
simplified relations).

The fact that $\theta_i^{(\cdot)}$ are all central enables a
rather quick simplification of the generators of $N$, by
translating the commutator relations into action by conjugation.
First,
\begin{eqnarray}\label{conj4}
\begin{array}{ccl}
{} 4_i \sp c_i \sp 4_i^{-1}  & = & {\theta_i^{(4,c)}} \sp c_i ,\\
{} 4_i \sp 7_i \sp 4_i^{-1}  & = & {\theta_i^{(4,7)}} \sp 7_i ,\\
{} 4_i \sp 10_i \sp 4_i^{-1} & = & \theta_i^{(4,10)} \sp c_i^3 \sp 10_i,\\
{} 4_i \sp 1_i \sp 4_i^{-1} & = & \theta_i^{(4,1)} \sp c^{-2}_i
\sp 7_i^{2} \sp 1_i, \end{array}
\end{eqnarray}
with trivial action of $4_i$ on every $\theta_i^{(\cdot)}$. Using
this we can remove $4_i$ from the relations; for example \eq{rel+}
becomes
\begin{eqnarray*}\label{rel++}
& & 1_i^{-1} c_i  7_i^{-1}  10_i^{-1} c_i^{-5} \sp 7_i^{2} 1_i 7_i
1_i^{-1} 7_i^{-2} c_i^5 10_i 7_i c_i^{-1}
1_i 7_i^{-1} \\
& & {} \quad\cdot \, 7_j 1_j^{-1} c_j
 7_j^{-1}  10_j^{-1}
c_j^{-5} 7_j^{2}  1_j 7_j^{-1}
 1_j^{-1} 7_j^{-2} c_j^5  10_j
7_j c_j^{-1}
1_j = e. %
\end{eqnarray*}
Again the first section is balanced on ${\theta_i^{(4,c)}},
{\theta_i^{(4,7)}}, {\theta_i^{(4,1)}}$ and ${\theta_i^{(4,10)}}$
so they do not show up in the simplified relation, and likewise
for the $\theta_{j}^{(\cdot)}$ in the second  section. This is
common among the $49$ generators of $N$, but there are some
exception. The effect of these exceptions is described later.

Next, we act by conjugation with $10_i$, applying the relations
\begin{eqnarray}\label{conj10}
\begin{array}{ccl}
{} 10_i \sp c_i \sp 10_i^{-1} & = & {\theta_i^{(10,c)}} \sp c_i\\
{} 10_i \sp 1_i \sp 10_i^{-1} & = & {\theta_i^{(10,1)}} \sp 1_i\\
{} 10_i \sp 7_i \sp 10_i^{-1} & = & {\theta_i^{(10,7)}} \sp 7_i,
\end{array}
\end{eqnarray} so continuing with the example, we obtain
\begin{eqnarray}\label{rel+++}
\begin{matrix}
& & {\theta_i^{(10,7)}}^{-1} 1_i^{-1} c_i  7_i^{-1}  c_i^{-5}
\sp 7_i^{2} 1_i 7_i 1_i^{-1} 7_i^{-2} c_i^5 7_i c_i^{-1}
1_i 7_i^{-1} \\
& & {} \quad\cdot \, {\theta_j^{(10,7)}} 7_j 1_j^{-1} c_j 7_j^{-1}
c_j^{-5} 7_j^{2}  1_j 7_j^{-1} 1_j^{-1} 7_j^{-2} c_j^5 7_j
c_j^{-1}
1_j = e.%
\end{matrix}
\end{eqnarray}
Next we act with $1_i$, using the relations
\begin{eqnarray}\label{conj1}
\begin{array}{ccl}
{}1_i \sp c_i \sp 1_i^{-1} & = & {\theta_i^{(1,c)}} \sp c_i\\
{}1_i \sp 7_i \sp 1_i^{-1}   & = & {\theta_i^{(1,7)}} \sp c_i^{-3}
\sp 7_i. \end{array}
\end{eqnarray}
To make the computation easier, we may at times conjugate the
relations (i.e. rotate each section), so the acting variables
envelope as shortest word as possible. Thus the relation
\eq{rel+++} is equivalent to
\begin{eqnarray*}
& & {\theta_i^{(10,7)}}^{-1} 1_i 7_i^{-1} 1_i^{-1} c_i  7_i^{-1}
c_i^{-5} \sp 7_i^{2} 1_i 7_i 1_i^{-1} 7_i^{-2} c_i^5 7_i c_i^{-1}
 \\
& & {} \quad\cdot \, %
{\theta_j^{(10,7)}} 1_j 7_j 1_j^{-1} c_j 7_j^{-1} c_j^{-5} 7_j^{2}
1_j 7_j^{-1} 1_j^{-1} 7_j^{-2} c_j^5 7_j c_j^{-1} = e,%
\end{eqnarray*}
which becomes
\begin{eqnarray*}
& & {\theta_i^{(10,7)}}^{-1} 7_i^{-1} c_i^{4} 7_i^{-1} c_i^{-5}
7_i^{2} c_i^{-3} 7_i^{-1} c_i^5 7_i c_i^{-1}
 \\
 & & {} \quad\cdot \, %
{\theta_j^{(10,7)}}  c_j^{-3} 7_j c_j 7_j^{-1} c_j^{-5} 7_j
c_j^{3}
7_j^{-2} c_j^5 7_j c_j^{-1} = e.%
\end{eqnarray*}
Finally, the last commutator relations translates to %
\begin{equation}\label{conj7}
7_i \sp c_i \sp 7_i^{-1} = {\theta_i^{(7,c)}} c_i, %
\end{equation}
so we obtain (among the $49$ relations)
\begin{eqnarray}\label{rel++++}
& & {\theta_i^{(10,7)}}^{-1} {\theta_i^{(7,c)}}
  {}\cdot \, %
{\theta_j^{(10,7)}} {\theta_j^{(7,c)}}^{-1} = e.%
\end{eqnarray}


It should be noted that some relations between the
$\theta_i^{(\cdot)}$ can be proved directly from their definition.
\begin{prop}\label{reltheta}
For every $i$, the elements $\theta_i^{(1,c)}, \theta_i^{(4,c)},
\theta_i^{(7,c)},\theta_i^{(10,c)}$ and $\theta_i^{(10,7)}$
satisfy
\begin{eqnarray*}
{} {\theta_i^{(10,c)}}& = & {\theta_i^{(10,7)}}^{-2}
{\theta_i^{(1,c)}}^{3},\\
{} {\theta_i^{(10,c)}}^{3} & = & e.\\
{} {\theta_i^{(7,c)}} & = & e,\\
{} {\theta_i^{(4,c)}}^{3} & = & e,\\
\end{eqnarray*}
\end{prop}
\begin{proof}
\forget These identities are obtained from acting with $4_i$ on
the relations \eq{conj10},\eq{conj1} and \eq{conj7}, from acting
with $10_i$ on \eq{conj1} and \eq{conj7}, and from acting with
$1_i$ on \eq{conj7}. %
\forgotten Conjugation of \eq{conj10}, \eq{conj1} and \eq{conj7}
by $4_i$ yields the four equations
\begin{eqnarray*}
{} {\theta_i^{(10,c)}}^{-2} {\theta_i^{(10,7)}}^{2}
 {\theta_i^{(1,c)}}^{-3}
{\theta_i^{(7,c)}}^{-6} & = &
e\\
{}{\theta_i^{(7,c)}}^{-3} & = & e,\\
{} {\theta_i^{(7,c)}}^{2} & = & e\\
{} {\theta_i^{(7,c)}}^{4}  & = & {\theta_i^{(4,c)}}^{3},
\end{eqnarray*}
and conjugation of $1_i\sp 7_i \sp 1_i^{-1} = {\theta_i^{(1,7)}}
\sp c_i^{-3} \sp 7_i$  (in \eq{conj1}) by $10_i$ yields
\begin{eqnarray*}
{} {\theta_i^{(10,c)}}^{3} & = & e.\\
\end{eqnarray*}
These five identities are equivalent to the ones listed in the
proposition.
\end{proof}

\section{The fundamental group of the Galois cover}

\subsection{A presentation for $K^*$}\label{sec:Kstar}

In the previous subsection we described the simplification process
of the $49$ relations defining $K^* = F^* / N$. Since we apply
relations that were already known modulo $\hat{N}$, it is not
surprising that the simplification removes all the original
generators $1_i,4_i,7_i,10_i$ and $c_i$, and leaves $49$ relations
on the $\theta_i^{(\cdot)} \in \hat{N}/N$ defined there.

\Eq{rel++++} is one example of these. This relation tells us that
${\theta_i^{(10,7)}} {\theta_i^{(7,c)}}^{-1}$ is independent of
$i$. Working out the relations in this manner, it turns out that
all the $15$ elements $\theta_i^{(\cdot)}$ are independent of the
lower index, namely $\theta_i^{(15)} = \theta_j^{(15)}$,
$\theta_i^{(17)} = \theta_j^{(17)}$, and so on. Therefore we may
remove the lower index and refer to the elements as
$\theta^{(15)}$, $\theta^{(17)}$ etc. As mentioned above, this was
already known when we started the computation, from Proposition
\ref{Nhat}.

Eventually, the $49$ relations in $N$ translate to the invariance
of all the $\theta_i^{(\cdot)}$, plus three more relations, which
refer to the four generators $\theta^{(1,c)}, \theta^{(7,c)},
\theta^{(10,c)},
\theta^{(10,7)}$: %
\begin{eqnarray*}
{\theta^{(10,c)}}^3 & = & e,\\
{\theta^{(1,c)}}^6 & = & e,\\
{\theta^{(10,c)}} & = & {\theta^{(10,7)}}^{-2}
{\theta^{(7,c)}}^{-8}
{\theta^{(1,c)}}^{-3}. %
\end{eqnarray*}

Combining this with Proposition~\ref{reltheta}, the relations on
the $\theta^{(\cdot)}$ can be summarized by
\begin{eqnarray}\label{thetarels}
\begin{array}{lcl}
{} {\theta^{(10,c)}}   & = & {\theta^{(10,7)}}^{4},\\
{} {\theta^{(10,7)}}^6 & = & {\theta^{(1,c)}}^{3},\\
{} {\theta^{(7,c)}}    & = & e,\\
{} {\theta^{(4,c)}}^3  & = & e,\\
{} {\theta^{(1,c)}}^6  & = & e.
\end{array}
\end{eqnarray}
It is worth noting that the five elements involved in the above
relations are all pure commutators (i.e.\ $\theta^{(1,c)} =
[1_i,c_i]$ etc., unlike $\theta^{(4,1)}$ or $\theta^{(1,7)} =
[1_i,7_i]c_i^3$ which are not commutators).


Let $\Omega' = \set{1,2,3,4,7,10,c,15,17,23}$ (i.e.\ $13$ is
replaced by $c$ in $\Omega$). Also let $\Omega_0 =
\set{c,7,1,10,4}$, ordered in this order (smaller first), and set
$\Theta = \set{\theta^{(\omega)} \suchthat \omega \in \Omega'
\minusset \Omega_0} \cup \set{\theta^{(\omega_1,\omega_2)}
\suchthat \omega_1,\omega_2 \in \Omega_0,\, \omega_1 < \omega_2}$,
the set of $15$ central elements discussed so far. We obtain the
following presentation:
\begin{cor}\label{Kstar}
The group $K^* = F^* / N$ is generated by the $18\cdot 10$
generators $\omega_i$ for $\omega \in \Omega'$ and $i =
1,\dots,18$ and the $15$ central generators $\Theta$, subject to
the following relation:
\begin{eqnarray*}
15_i & = & \theta^{(15)} \sp 7_i^{-1} \sp 4_i, \\ 
17_i & = & \theta^{(17)} \sp 10_i \sp7_i, \\ 
23_i & = & \theta^{(23)} \sp c_i \sp 7_i^{-1}, \\ %
2_i & = & \theta^{(2)} \sp c_i^2 \sp 7_i^{-1} \sp 10_i^{-1}\sp 1_i, \\ 
3_i & = & \theta^{(3)} \sp c_i \sp 7_i^{-1}\sp 10_i^{-1} \sp 1_i,\\
{}[7_i,c_i] & = & \theta^{(7,c)}, \\
{}[1_i,c_i] & = & \theta^{(1,c)}, \\
{}[1_i,7_i] & = & \theta^{(1,7)} \sp c_i^{-3}, \\
{}[10_i,c_i] & = & \theta^{(10,c)}, \\
{}[10_i,7_i] & = & \theta^{(10,7)}, \\
{}[10_i,1_i] & = & \theta^{(10,1)}, \\
{}[4_i,c_i] & = & \theta^{(4,c)}, \\
{}[4_i,1_i] & = & \theta^{(4,1)} \sp c_i^{-2}\sp 7_i^2, \\
{}[4_i,7_i] & = & \theta^{(4,7)}, \\
{}[4_i,10_i] & = & \theta^{(4,10)} \sp c_i^{3},
\end{eqnarray*}
and
\begin{eqnarray}
{}\theta^{(10,c)} & = & {\theta^{(10,7)}}^{4} \nonumber, \\ %
{\theta^{(10,7)}}^{6} & = & {\theta^{(1,c)}}^{3}, \nonumber\\
{}{\theta^{(7,c)}}    & = & e, \label{thetaincor}\\
{}{\theta^{(4,c)}}^3  & = & e, \nonumber \\
{}{\theta^{(1,c)}}^6  & = & e. \nonumber
\end{eqnarray}
\forget $15$ expressions being independent of $i$ and central:
\begin{eqnarray*}
15_i \sp (7_i^{-1} \sp 4_i)^{-1}, \quad 17_i \sp(10_i
\sp7_i)^{-1}, \quad 23_i \sp (c_i \sp
7_i^{-1})^{-1},  \\
2_i \sp (c_i^2 \sp 7_i^{-1} \sp 10_i^{-1} \sp 1_i)^{-1}, \quad 3_i
\sp (c_i \sp 7_i^{-1} \sp 10_i^{-1} \sp 1_i)^{-1}, \\
{}[c_i,1_i], \quad [c_i,7_i], \quad [c_i,4_i], \quad [c_i,10_i], \\
{} [1_i,10_i], \quad  [10_i,7_i], \quad [7_i,4_i], \\
{}[7_i,1_i] \sp c_i^{-3}, \quad [4_i,10_i] \sp c_i^{-3}, \\
{}[4_i,1_i] \sp (c_i^{-2}\sp 7_i^2)^{-1},
\end{eqnarray*}
plus the following five relations:
\begin{eqnarray*}
{}[c_i,10_i]   & = & [10_i,7_i]^{2} [c_i,1_i]^3. \\ %
{}[c_i,1_i]^6  & = & e, \\
{}[c_i,4_i]^3  & = & e, \\
{}[c_i,7_i]    & = & e, \\
{}[c_i,10_i]^3 & = & e.
\end{eqnarray*}
\forgotten
\end{cor}

\begin{prop}
The subgroup $\sg{\Theta}$ of $K^*$ is isomorphic to $\Z^{10}
\times (\Z/3)^2 \times \Z/12$.
\end{prop}
\begin{proof}
It is easy to see that the \emph{abstract} group generated by
$\Theta$, subject only to the relations given in \Eq{thetaincor},
is indeed isomorphic to $\Z^{10} \times (\Z/3)^2 \times \Z/12$.

In order to show that the natural map from the abstract group
$\sg{\Theta}$ into $K^*$ is indeed an embedding, we need to
construct $K^*$ as an extension of $\sg{\Theta}$. This can easily
be done as a sequence of HNN extensions. First adjoin the elements
$c_i$ (all inducing the trivial automorphism on $\sg{\Theta}$).
Then adjoin $7_1,\dots,7_{18}$ with $7_i$ commuting with every
generator $\omega_j$ ($j\neq i$) and with the elements of
$\Theta$, and acting on $c_i$ by $c_i \mapsto \theta^{(7,c)} c_i$.
The automorphisms commute, and so we can form the (repeated) HNN
extension by $7_1,\dots,7_{18}$.

In a similar manner we adjoin the $1_i$, then the $10_i$ and
finally the $4_i$. In each step one needs to verify that the  map
we want to be induced is indeed an automorphism, but this was done
(implicitly) in the proof Proposition \ref{reltheta}.
\end{proof}

\begin{rem}
The quotient $K^*/\sg{\Theta}$ is isomorphic to $H^{18}$ (by
inspection, see the definition of $H$).
\end{rem}

\begin{prop}\label{CentofKstar}
The center of $K^*$ is generated by $\Theta$ and
$c_1^6,\dots,c_{18}^6$. In particular $\Cent(K^*) \isom \Z^{18}
\times \sg{\Theta}$, and $K^*/\Cent(K^*) \isom (H/\sg{c^6})^{18}$.
\end{prop}
\begin{proof}
Since $K^*$ can be viewed as an iterated HNN extension of
$\sg{\Theta}$ by the $c_i, 7_i, 1_i, 10_i$ and $4_i$, every
element can be written as a product $\theta_0
\prod_{i=1}^{18}{c_i^{n_{c,i}} 7_i^{n_{7,i}} 1_i^{n_{1,i}}
10_i^{n_{10,i}} 4_i^{n_{4,i}}}$ for $\theta_0 \in \sg{\Theta}$.
Assume such an element is central. Conjugating by $7_i$ proves
that $n_{1,i} = n_{4,i} = 0$, since $\theta^{(1,7)}$ and
$\theta^{(4,7)}$ have infinite order. Likewise conjugation by
$4_i$ shows $n_{7,i} = n_{10,i} = 0$. Eventually conjugation by
$1_i$ shows $n_{c,i}$ is divisible by $6$, and on the other hand
it is easy to see that $c_i^6$ are central.
\end{proof}

Since $(F^*/N) / (\hat{N}/N) = F^*/\hat{N}$ is known to be
isomorphic to $H^{18}$ (Subsection \ref{groupH}), we have a short
exact sequence
$$1 \llra \Z^{10}
\times (\Z/3)^2 \times \Z/12 \llra K^* \llra H^{18} \llra 1.$$ In
particular from Corollary \ref{Hisnilp} we conclude that $K^*$ is
nilpotent of class at most $4$. However, we have

\begin{prop}\label{Kstarnilp}
$K^*$ is nilpotent of class $3$. In more details: \\
a. The commutator subgroup of $K^*$ is
$$[K^*,K^*] = \sg{\theta^{(1,c)},
\theta^{(1,7)} \sp c_i^{-3},
 \theta^{(10,7)}, \theta^{(10,1)}, \theta^{(4,c)}, \theta^{(4,1)} \sp c_i^{-2}\sp 7_i^2,
\theta^{(4,7)}, \theta^{(4,10)} \sp c_i^{3}},$$ isomorphic to
$\Z^5 \times (\Z/3)^2 \times \Z/12$.
\\
b. The next term in the upper central series is
$$[K^*,[K^*,K^*]]
= \sg{{\theta^{(1,c)}}^3, {\theta^{(1,c)}}{\theta^{(1,7)}}^2
c_i^{-6} , {\theta^{(4,c)}}^{}{\theta^{(4,7)}}^{2}},$$ isomorphic
to $\Z^2 \times \Z/2$.
\\
c. $[K^*,[K^*,[K^*,K^*]]] = 1$ and in particular
$[[K^*,K^*],[K^*,K^*]] = 1$.
\\
d. $[K^*,K^*] / [K^*,[K^*,K^*]] \isom \Z^3 \times (\Z/6)^3$.
\end{prop}
\begin{proof}
Direct computation from the defining relations of $K^*$.
\end{proof}

\begin{rem}
The abelianization of $K^*$ is isomorphic to $\Z^{61}\times
(\Z/6)^{17}$.
\end{rem}
\begin{proof}
The abelianization is generated by $\Theta$ and the $5\cdot 18$
elements $c_i,1_i,7_i,10_i,4_i$. Setting the commutators to be
trivial, we are left only with $\theta^{(2)},\theta^{(3)},
\theta^{(15)}, \theta^{(17)}, \theta^{(23)}$ and $\theta^{(1,7)},
\theta^{(4,1)}, \theta^{(4,10)}$ from $\Theta$, with the relations
${\theta^{(4,10)}}^{-1} \equiv \theta^{(1,7)} \equiv c_i^3$ and
$c_i^2 \equiv \theta^{(4,1)}7_i^2$. {}This is equivalent to $c_i
\equiv \theta^{(1,7)} {\theta^{(4,1)}}^{-2} 7_i^{-2}$ and %
\begin{equation}\label{abel7}
7_i^6 \equiv {\theta^{(1,7)}}^{2}{\theta^{(4,1)}}^{-3}.
\end{equation}

This can further be replaced by $(7_i7_1^{-1})^6 \equiv 1$ for $i
> 1$, and ${\theta^{(4,1)}} \equiv (7_1^{-3}
{\theta^{(1,7)}}{\theta^{(4,1)}}^{})^{-1}$, so by a base change we
obtain the independent generators $\theta^{(2)},\theta^{(3)},
\theta^{(15)}, \theta^{(17)}, \theta^{(23)}$, $1_i,4_i,10_i$
($i=1,\dots,18$), $7_1$, $7_1^{-3}
{\theta^{(1,7)}}^{}{\theta^{(4,1)}}^{}$, and the $17$ elements
$7_i7_1^{-1}$ ($i > 1$). The result follows.
\end{proof}

\subsection{Back to $K$}\label{back}

Recall from Subsection \ref{modcent} that $\pitil^{\Aff} \isom S_n
\semidirect K$ where $K$ is the kernel of the map $K^* \ra
\Z^{10}$ defined as the restriction of the map from $F^*$ (see
Subsection \ref{ss:struc}).

\begin{cor}\label{2.10}
$K$ is a normal subgroup of the group $K^*$, whose presentation is
given in Corollary \ref{Kstar}. In fact $K$ is the kernel of the
map $\ab \co K^* \ra \Z^{10} = \sg{e_2,e_3,\dots,e_{23}}$ induced
from $F^* \ra \Z^{10}$.

Explicitly, the map $\ab$ is defined by sending the generators
$\omega_i$ to $e_\omega$ for every $\omega \in \Omega'$, with
\begin{eqnarray*}
\theta^{(15)} & \mapsto & e_{15} + e_{7} - e_{4}\\
\theta^{(17)} & \mapsto & e_{17} - e_{7} - e_{10}\\
\theta^{(23)} & \mapsto & e_{23} - e_{c} + e_{7}\\
\theta^{(2)} & \mapsto & e_{2} - 2 e_{c} + e_{7} + e_{10} - e_{1}\\
\theta^{(3)} & \mapsto & e_{2} - e_{c} + e_{7} + e_{10} - e_{1}\\
\theta^{(1,7)} & \mapsto & 3 e_{c}\\
\theta^{(4,1)} & \mapsto & 2 e_{c} - 2 e_7\\
\theta^{(4,10)} & \mapsto & - 3 e_{c},
\end{eqnarray*}
and the other generators in $\Theta$ map to $0$.
\end{cor}

\begin{rem}\label{2.11}
The commutator subgroup of $K$ is $[K,K] = [K^*,K^*]$. It follows
that $[K^*,[K^*,K^*]] = [K^*,[K,K]] \,\sub\, [K,[K^*,K^*]] =
[K,[K,K]]$. Therefore $K$ is nilpotent of class $3$ (and not
less).
\end{rem}
\begin{proof}
To prove that $[K^*,K^*] \sub [K,K]$, let $\omega_i,\omega_i'$ be
any two generators of $K^*$. Let $j,k \neq i$ be two distinct
indices, then $\omega_i\omega_j^{-1}, \omega'_{i} {\omega'_k}^{-1}
\in K$ and $[\omega_i\omega_j^{-1}, \omega'_{i} {\omega'_k}^{-1}]
= [\omega_i,\omega'_i]$.

For the corollary, note that $[B,[A,A]] \sub [A,[B,B]]$ for every
two groups $A \sub B$ by the three subgroups lemma.
\end{proof}

Since $\pitil^{\Aff} \isom S_n \semidirect K$, we proved our main
result:
\begin{thm}
$\pitil^{\Aff}$ is virtually nilpotent of class $3$.
\end{thm}

\begin{prop}\label{2.13}
The centralizer of $K$ in $K^*$ is generated by $\Cent(K^*)$ and
the element $c_1c_2\cdots c_{18}$.
\end{prop}
\begin{proof}
We know that $K$ is generated by the elements
$\omega_\alpha\omega_\beta^{-1}$ for $\omega \in \Omega$ and
$\alpha,\beta = 1,\dots,18$. Mimicking the argument of Proposition
\ref{CentofKstar}, we conclude that elements of the centralizer
are all of the form $$z \cdot (\prod{c_i})^{n_c} (\prod 7_i)^{n_7}
(\prod 1_i)^{n_1} (\prod 10_i)^{n_{10}} (\prod 4_i)^{n_4}$$ for $z
\in \Cent(K^*)$. Conjugating by $7_\alpha 7_{\beta}^{-1}$, we
conclude that $n_{1} = 0$. Likewise conjugation by $10_\alpha
10_{\beta}^{-1}, 1_\alpha 1_{\beta}^{-1}$ and $4_\alpha
4_{\beta}^{-1}$ shows that $n_{4} = 0$, $n_{10} = 0$ and $n_{7} =
0$. On the other hand, $\prod c_i$ commutes with $K$.
\end{proof}

\ifXY
\begin{figure}
\begin{equation*}
\xymatrix@C=12pt@R=10pt{ 
    {}
        &
        & K^* \ar@{-}[d] 
\\
    {}
        & K \Cent(K^{*}) \ar@{-}[dd] \ar@{-}[ld]
        & K \op{C}_{K^*}(K) \ar@{-}[dd] \ar@{=}[l]
\\
    K \ar@{-}[dd]
        &
        &
\\
    {}
        & \Cent(K) \Cent(K^*) \ar@{-}[dd] \ar@{-}[ld]
        & \op{C}_{K^*}(K) \ar@{=}[l]
\\
    \Cent(K) \ar@{-}[dd]
        &
        &
\\
    {}
        & \Cent(K^*) \ar@{-}[ld] 
        & 
\\
    K \cap \Cent(K^*)
        &
        &
}
\end{equation*}
\caption{Subgroups of $K^*$ and $K$}\label{KK}
\end{figure}
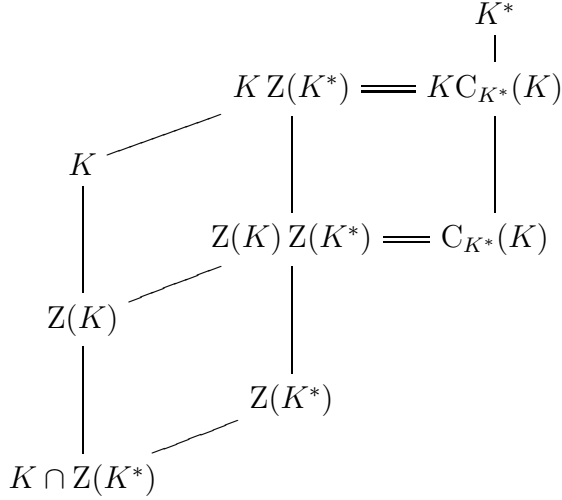
\else
\fi 

\begin{cor}
The centralizer of $K$ in $K^*$ is $\operatorname{C}_{K^*}(K) =
\Cent(K) \Cent(K^*)$.
\end{cor}
\begin{proof}
The inclusion $\supseteq$ is trivial, and by the last proposition
it remains to prove that $\hat{c} = c_1 \dots c_{18} \in
\Cent(K)\Cent(K^*)$. But $\ab(\hat{c} {\theta^{(4,10)}}^6) = 0$,
so  $\hat{c} {\theta^{(4,10)}}^6 \in K \cap \Cent(K^*) \sub
\Cent(K)$ and $\hat{c} = (\hat{c} {\theta^{(4,10)}}^6)
{\theta^{(4,10)}}^{-6} \in \Cent(K)\Cent(K^*)$.
\end{proof}

\subsection{The projective relation}

\ifXY
\begin{figure}
\begin{equation*}
\xymatrix@C=24pt@R=24pt{ 
    \pi_1(\GalAff{X}) \ar@{^(->}[r]
    \ar@{->}[d]
        & \pitil^{\Aff} = \pi_1(\C^2 - S) / \sg{\Gamma_j^2,\Gamma_{j'}^2}
        \ar@{->}[r] \ar@{->}[d]
        & S_n \ar@{=}[d]
\\
    \pi_1(\Gal{X}) \ar@{^(->}[r]
        & \pitil = \pi_1(\C\P^2 - S) / \sg{\Gamma_j^2,\Gamma_{j'}^2}
        \ar@{->}[r]
        & S_n
}
\end{equation*}
\caption{The affine and projective parts}\label{proj}
\end{figure}
\else
\fi 

We denote the affine part of the Galois cover of $\Gal{X}$ by
$\GalAff{X}$. Recall from the introduction that the fundamental
group of  $\Gal{X}$ is the kernel of the projection $\pitil =
\pi_1(\C\P^2 - S)/\sg{\Gamma_j^2,\Gamma_{j'}^2} \ra S_n$, namely
the bottom line of the diagram in Figure \ref{proj} is exact. In a
similar manner, the upper line is also exact, and so we know from
\Eq{pitilis} that $\pi_1(\GalAff{X}) = K = F/N$ (where $F$ and $N$
are defined in Subsection \ref{FN}).

Let $p$ denote the `projective product' $p = \Gamma_{27'}
\Gamma_{27} \dots \Gamma_{1'} \Gamma_1$ (see \Eq{projrel}). {}The
general theory guarantees that
$$1 \llra \sg{p} \llra \pitil^{\Aff} \llra \pitil \llra 1$$ is a short exact
sequence (and in fact a central extension, see
\cite[Prop.~3.2]{Te7}). 

Since $p \in K^*$ is clearly mapped to the trivial permutation in
$S_n$, it belongs to the kernel $K$, and so the following is again
a central extension:
$$1 \llra \sg{p} \llra K \llra \pi_1(\Gal{X}) \llra 1,$$
so $\pi_1(\Gal{X})$ (which is the group we are really after) is
isomorphic to $K / \sg{p}$.
\begin{thm}\label{main}
The fundamental group $\pi_1(\Gal{X})$ is nilpotent,  of
nilpotency class $3$.
\end{thm}
\begin{proof}
In Remark \ref{2.11} we saw that $K$ is nilpotent of class $3$,
and that $[K,[K,K]] \isom \Z^2 \times \Z/2$ (by Proposition
\ref{Kstarnilp}.b). Therefore $$[K/\sg{p},[K/\sg{p},K/\sg{p}]] =
[K,[K,K]]\sg{p}/\sg{p}$$ cannot be trivial.
\end{proof}

In order to compute $p$ explicitly, we first apply the
substitutions given in Subsection \ref{remove}, obtaining a word
of length $3822$ on the $27$ generators $\Delta$ of \Eq{Delta}.
This word can be viewed as an element of $\CoxY{T}$ (modulo the
relations defining $\pitil$). Next, since $p$ maps to the trivial
permutation, it is an element of $F^{*}/N$, namely a word in the
generators $\omega_{\alpha}$ ($\omega \in \Omega$ and $\alpha =
1,\dots,18$); in fact $p \in K = F/N$ since $p$ comes from
$\CoxY{T}$. However, $F/N$ is a central extension of $F/\hat{N}
\isom H^{18}$ by $\hat{N}/N$ (see Proposition \ref{Nhat}), so we
can decompose $p$ as a product of elements from these groups. Note
that $p$ being central in $\pitil$ guarantees that $p$ is
invariant under the action of $S_n$, and so it has to have the
same component (modulo $\hat{N}$) in each factor of the $H_1
\times \cdots \times H_{18}$.

The substitutions made so far yield an element $p$ of length
$325$, as follows:
\def\spe{{\,}}
\begin{eqnarray*}
p & = & {10}_{1}\spe 7_{1}\spe 1^{-1}_{1}\spe 17_{1}\spe 23_{1}
2_{1}\spe {4}^{-1}_{1}\spe 7^{-1}_{1}\spe 3^{-1}_{1}\spe
17^{-1}_{1}\spe 13^{-1}_{1}\spe 2_{1}
\\\spe %
& & \spe \cdot\spe 17_{2}\spe 3_{2}\spe 1^{-1}_{2}\spe \\\spe %
& & \spe \cdot\spe 17_{3}\spe 3_{3}\spe 1^{-1}_{3}\spe \\\spe %
& & \spe \cdot\spe {10}^{-1}_{4}\spe 2^{-1}_{4}\spe
23^{-1}_{4}\spe 17^{-1}_{4}\spe 1_{4}\spe 7^{-1}_{4}\spe
3^{-1}_{4}\spe 17^{-1}_{4}\spe 1_{4}\spe {4}_{4}\spe 2^{-1}_{4}
13_{4}
1_{4}\spe 7_{4}\spe 1^{-1}_{4}\spe 17_{4}\spe 3_{4}\spe \\\spe %
& & \spe \cdot\spe {4}^{-1}_{5}\spe 15_{5}\spe 2^{-1}_{5}
13_{5}\spe 17_{5}\spe 3_{5}\spe 7_{5}
1^{-1}_{5}\spe 13^{-1}_{5}\spe 2_{5}\spe \\\spe %
& & \spe \cdot\spe {4}^{-1}_{6}\spe 1^{-1}_{6}\spe 17_{6}
3_{6}\spe 7_{6}\spe {4}_{6}
2^{-1}_{6}\spe 13_{6}\spe 1_{6}\spe 7^{-1}_{6}\spe 1^{-1}_{6}\spe 13^{-1}_{6}\spe 2_{6}\spe \\\spe %
& & \spe \cdot\spe {4}^{-1}_{7}\spe 1^{-1}_{7}\spe 13^{-1}_{7}
15^{-1}_{7}\spe 1_{7}\spe {4}_{7}\spe 2^{-1}_{7}\spe 13_{7}\spe
1_{7}\spe 7^{-1}_{7}\spe 3^{-1}_{7}\spe 17^{-1}_{7}
13^{-1}_{7}\spe 2_{7}\spe \\\spe %
& & \spe ...\spe \\
\forget\spe & & \spe {10}_{8}\spe 7^{-1}_{8}\spe 13_{8}\spe
1_{8}\spe 7^{-1}_{8}\spe 3^{-1}_{8}\spe 17^{-1}_{8}\spe
13^{-1}_{8}\spe 2_{8}\spe {4}^{-1}_{8}
1^{-1}_{8}\spe 13^{-1}_{8}\spe 15^{-1}_{8}\spe {4}_{8}\spe \\\spe %
& & \spe {10}^{-1}_{9}\spe 2^{-1}_{9}\spe 23^{-1}_{9}
17^{-1}_{9}\spe 1_{9}\spe {4}_{9}\spe 2^{-1}_{9}
13_{9}\spe 1_{9}\spe \\\spe %
& & \spe {4}^{-1}_{10}\spe 7^{-1}_{10}\spe 7^{-1}_{10}\spe
15_{10}\spe 2^{-1}_{10}\spe 13_{10}\spe 17_{10}\spe 3_{10}\spe
7_{10}\spe 7_{10}\spe 7_{10}
1^{-1}_{10}\spe 13^{-1}_{10}\spe 2_{10}\spe \\\spe %
& & \spe {4}_{11}\spe 2^{-1}_{11}\spe 13_{11}\spe 1_{11}
7_{11}\spe 1^{-1}_{11}\spe 7^{-1}_{11}\spe 3^{-1}_{11}
17^{-1}_{11}\spe 13^{-1}_{11}\spe 2_{11}\spe 15^{-1}_{11}\spe
1_{11}\spe 7^{-1}_{11}\spe 1^{-1}_{11}\spe 13^{-1}_{11}\spe 2_{11}
{4}^{-1}_{11}\spe 1^{-1}_{11}\spe 15_{11}\spe 13_{11}\spe 1_{11}
{4}_{11}\spe 2^{-1}_{11}\spe 13_{11}\spe 17_{11}\spe 3_{11} 7_{11}
\\\spe %
& & \spe {4}^{-1}_{12}\spe 1^{-1}_{12}\spe 7^{-1}_{12}
3^{-1}_{12}\spe 17^{-1}_{12}\spe 13^{-1}_{12}\spe 2_{12}\spe
13_{12}\spe 1_{12}\spe {4}_{12}\spe 2^{-1}_{12}\spe 13_{12}
17_{12}\spe 3_{12}\spe 7_{12}\spe 1^{-1}_{12}\spe 13^{-1}_{12}
2^{-1}_{12}\spe 13_{12}\spe 17_{12}\spe 3_{12}\spe 7_{12}\spe
1_{12}\spe 7^{-1}_{12}
1^{-1}_{12}\spe 13^{-1}_{12}\spe 2_{12}\spe \\\spe \spe %
& & \spe {10}_{13}\spe 7_{13}\spe 1^{-1}_{13}\spe 17_{13}
3_{13}\spe 7_{13}\spe 1^{-1}_{13}\spe 17_{13}\spe 23_{13}\spe
13_{13}\spe 1_{13}\spe 7^{-1}_{13}\spe 3^{-1}_{13}\spe
17^{-1}_{13}\spe 13^{-1}_{13}\spe 2_{13}\spe {4}^{-1}_{13}
7^{-1}_{13}\spe 3^{-1}_{13}\spe 17^{-1}_{13}
13^{-1}_{13}\spe 2_{13}\spe \\\spe %
& & \spe {10}^{-1}_{14}\spe {4}^{-1}_{14}\spe 15_{14}
{10}_{14}\spe {4}_{14}\spe 2^{-1}_{14}\spe 13_{14}\spe 1_{14}
7_{14}\spe 1^{-1}_{14}
13^{-1}_{14}\spe 2_{14}\spe {4}^{-1}_{14}\spe 1^{-1}_{14}\spe 17_{14}\spe 3_{14}\spe \\\spe %
& & \spe {10}_{15}\spe 7_{15}\spe 1^{-1}_{15}\spe 17_{15}
3_{15}\spe 7_{15}\spe 1^{-1}_{15}\spe 17_{15}\spe 7_{15}
1^{-1}_{15}\spe 13^{-1}_{15}\spe 2_{15}\spe {4}^{-1}_{15}
1^{-1}_{15}\spe 13^{-1}_{15}\spe 2_{15}\spe {4}^{-1}_{15}
1^{-1}_{15}\spe 17_{15}\spe 3_{15}\spe {10}^{-1}_{15}
{4}^{-1}_{15}\spe 15_{15}\spe {10}_{15}\spe 7^{-1}_{15}\spe
13_{15}\spe 1_{15}\spe 7^{-1}_{15}\spe 3^{-1}_{15}\spe
17^{-1}_{15}\spe 1_{15}\spe 7^{-1}_{15}\spe 3^{-1}_{15}
17^{-1}_{15}\spe 1_{15}
{4}_{15}\spe \\\spe %
& & \spe {10}^{-1}_{16}\spe 15^{-1}_{16}\spe {4}_{16}
{10}_{16}\spe 3^{-1}_{16}\spe 17^{-1}_{16}\spe 1_{16}\spe
{4}_{16}\spe 2^{-1}_{16}\spe 13_{16}\spe 1_{16}\spe {4}_{16}
2^{-1}_{16}\spe 13_{16}\spe 1_{16}\spe 7^{-1}_{16}\spe
17^{-1}_{16}\spe 1_{16}\spe 7^{-1}_{16}\spe 3^{-1}_{16}
17^{-1}_{16}\spe 1_{16}\spe 7^{-1}_{16}\spe {10}^{-1}_{16}
{4}^{-1}_{16}\spe 1^{-1}_{16}\spe 17_{16}\spe 3_{16}\spe 7_{16}
1^{-1}_{16}\spe 13^{-1}_{16}\spe 7_{16}\spe {10}^{-1}_{16}
15^{-1}_{16}\spe {4}_{16}\spe {10}_{16}\spe 3^{-1}_{16}
17^{-1}_{16}\spe 1_{16}
{4}_{16}\spe 2^{-1}_{16}\spe 13_{16}\spe 1_{16}\spe {10}^{-1}_{16}\spe 7_{16}\spe \\\spe %
& & \spe {10}^{-1}_{17}\spe 2^{-1}_{17}\spe 13_{17}\spe
17_{17}\spe 3_{17}\spe 7_{17}\spe 1^{-1}_{17}\spe 13^{-1}_{17}
2_{17}\spe {10}_{17}\spe 3^{-1}_{17}\spe 1_{17}\spe 7^{-1}_{17}
1^{-1}_{17}\spe 13^{-1}_{17}\spe 2_{17}\spe {4}^{-1}_{17}
7^{-1}_{17}\spe 3^{-1}_{17}\spe 17^{-1}_{17}\spe 1_{17}\spe 7^{-1}_{17}\spe \\\spe %
\forgotten\spe & & \spe \cdot\spe {4}_{18}\spe 2^{-1}_{18}
13_{18}\spe 1_{18}\spe 7_{18}\spe 1^{-1}_{18}\spe 7_{18}
1^{-1}_{18}\spe 13^{-1}_{18}\spe 2_{18}\spe {4}^{-1}_{18}\spe
7_{18}\spe {4}_{18}\spe 2^{-1}_{18}\spe 13_{18}\spe 1_{18}
7^{-1}_{18}\spe 17^{-1}_{18}\spe 1_{18} %
\end{eqnarray*}
We then compute each component using the relations of $K^*$ given
in Subsection \ref{sec:Kstar}, and find that
\begin{equation}\label{pis}
p = {\theta^{(4,10)}}^3 {\theta^{(4,c)}}^2 {\theta^{(1,7)}}^{-3}
{\theta^{(1,c)}} c_1 c_2 \dots c_{18}.%
\end{equation} %
As expected, $p$ is invariant under the action of $S_n$,
centralizes $K$ by Proposition \ref{2.13}, and $\ab(p) = 0$ by
Corollary \ref{2.10}, showing that $p \in K$.

\forget 
Since $K^*/\Cent(K^*) \isom H^{18}$, we have that $K^* /
\sg{p}\Cent(K^*) \isom H_0^{18}$ where $H_0 = H/ \sg{c}$.
Following Subsection \ref{groupH}, $H_0$ is equal to $\sg{1,4,7,10
\subjectto [4,1] = 7^2}$ with all the other pairs of generators
commute. In particular $H_0$ is a central extension of $\Z^3$ by
$\Z = \sg{7}$. Moreover, $$H_0 / H_0' = \sg{1,4,7,10 \subjectto
7^2 = e} \isom \Z^3 \times \Z/2.$$ \forgotten

\forget 
Note that since $p \in K$, the map $\ab$ discussed in
Subsection \ref{back} takes $p$ to $0$, and so it is still well
defined on $K^* / \sg{p}\Cent(K^*)$. We thus have the
following version of Corollary \ref{KmodCent}: %
\begin{cor}\label{KmodCentp}
$K / \sg{p}(K \cap \Cent(K^*))$ is the kernel of the map from
$K^*/\sg{p}\Cent{K^*}$ to $\Z^3 \times \Z/2$, defined as before.
\end{cor}
\forgotten

\forget 

\section{$\pi_1(\Gal{X})$ modulo $\Gamma_j = \Gamma_{j'}$}

Recall that $\pi_1(\GalAff{X})$ is generated by the
$\Gamma_j,\Gamma_{j'}$, and can be identified with the group $K$
defined in Subsection \ref{modcent}. There is a group $K^*$
containing $K$, with $K^*/K \isom \Z^{10}$, and $K^*/\Cent(K^*)
\isom H^{18}$ ($H$ is defined in Subsection \ref{groupH}).

As mentioned in the introduction, in a previous paper we computed
the quotient $C$ of $\pi_1(\Gal{X})$ defined by identifying
$\Gamma_{j} = \Gamma_{j'}$ for every $j$. In this section we
compute $C$ in the new context; this can serve as another way to
verify the computations.

To find the defining relations of $C$ in terms of the generators
of $F^*$, we perform (what is by now) the usual algorithm. The
$27$ relations $\Gamma_j = \Gamma_{j'}$ are transformed using the
substitutions of Subsection \ref{remove} to $27$ words on the
generators $\Delta$ of \Eq{Delta}, of total length $3824$.

\begin{exmpl}
In Subsection \ref{remove} we report that $8' = 6\sp 7\sp 12\sp
20\sp 24\sp 20\sp 12\sp 7\sp 6$. Therefore, in $C$ we have the
relation $6\sp 7\sp 12\sp 20\sp 24\sp 20\sp 12\sp 7\sp 6 = 8$.

Applying the map $\Phi$ of \Eq{pitilis}, this relation becomes
$$(2\sp 3) (2\sp 6) 7_6^{-1} 7_2 (6 \sp 14) (10 \sp 18) (14 \sp
18)(10 \sp 18) (6 \sp 14) (2\sp 6) 7_6^{-1} 7_2 (2\sp 3) = (3 \sp
10),$$ since all the generators here except for $\Gamma_7$ belong
to the spanning subtree $T_0$, and $7$ points from $\alpha_7 = 2$
to $\beta_7 = 6$ in Figure \ref{subtree}.

Simplifying this relation via the action of $S_{18}$ on $F^{*}$
(and hence on $K^{*} = F^{*}/N$), we obtain the relation $7_3 =
7_{10}$, which implies $$7_i = 7_j$$ for every $i \neq j$ (in the
group $C$, of course).
\end{exmpl}

The outcome of this computation, carried out for all the $27$
relations defining $C$, is that $C$ is isomorphic to the quotient
of $F/N \sub F^*/N = \sg{\omega_\alpha}_{\omega\in \Omega, \alpha
= 1,\dots,18}/N$, with respect to the following elements being
invariant under the symmetric group (and therefore do not require
subscripts): $7, 10, 17, 23, 15\sp 4^{-1}, 3\sp 1^{-1}, 2\sp
1^{-1}, 13\sp 1 \sp 4 \sp 1^{-1}, 1 \sp 4 \sp 1^{-1} \sp 4^{-1}$.
As before, the definition of $F^*$ implies that (`short')
invariant elements are in particular central.

((I computed the $49$ relations which generate $N$ -- and modulo
the above list of central element, they all become trivial. Thus
$N \sub U$ where $U$ is the normal subgroup generated by $j j'$.))

\forgotten

\section*{Acknowledgments}

The first named author is partially supported by
 EU-network HPRN-CT-2009-00099 (EAGER), the Emmy Noether Research
 Institute for Mathematics and the Minerva Foundation of Germany,
and an Israel Science Foundation grant \#8008/02-3 (Excellency
Center ``Group Theoretic Methods in the Study of Algebraic
Varieties").
 The work was also supported in part by the Edmund Landau Center
 for Research in Mathematical Analysis and Related Areas, sponsored
 by the Minerva Foundation (Germany). The first named author
 thanks the Landau center and especially
 Prof. Ruth Neumark-Lawrence for their hospitality.

\section*{Appendix}

The text refers to four files, each containing the presentation of
the fundamental group or one of its sections. The files
\begin{itemize}
\item[{[App1]:}] \software{presentation.txt}
\item[{[App2]:}] \software{presentation27.txt}
\item[{[App3]:}] \software{relations\_of\_N.txt}
\item[{[App4]:}] \software{relations\_of\_H.txt}
\end{itemize}
can be downloaded from
\software{http://www.math.biu.ac.il/\~{}vishne/downloads/TxT/}.

\end{document}